\documentclass[10pt]{amsart}
\usepackage[T1]{fontenc}
\usepackage[margin=1.5in]{geometry}
\usepackage[utf8]{inputenc}
\usepackage{xcolor}
\usepackage{amstext}
\usepackage{amsthm}
\usepackage{amssymb}
\usepackage{stmaryrd}
\usepackage{geometry}
\usepackage[bookmarks=false,
 breaklinks=false,pdfborder={0 0 1},backref=false,colorlinks=true]
 {hyperref}
\hypersetup{
 linkcolor=blue,citecolor=blue,urlcolor=blue}

\makeatletter
\numberwithin{equation}{section}
\numberwithin{figure}{section}

\@ifundefined{date}{}{\date{}}
\usepackage{amsthm}\numberwithin{equation}{section}
\newtheorem{theorem}{Theorem}[section]
\newtheorem{proposition}[theorem]{Proposition}\newtheorem{lemma}[theorem]{Lemma}\newtheorem{corollary}[theorem]{Corollary}\theoremstyle{definition}
\newtheorem*{problem}{Problem}\theoremstyle{remark}
\newtheorem{remark}[theorem]{Remark}\newcommand{\Scal}{\operatorname{Scal}}
\newcommand{\Ric}{\operatorname{Ric}}
\newcommand{\tr}{\operatorname{tr}}
\newcommand{\Id}{\operatorname{Id}}
\newcommand{\CP}{\mathbb{CP}}
\newcommand{\R}{\mathbb{R}}
\newcommand{\Q}{\mathbb{Q}}
\title[Yau's pinching problem]{Yau's pinching problem and Vanishing of Betti Numbers}
\author{Xiaolong Li}
\address{Department of Mathematics and Statistics, Auburn University,
Auburn, AL 36849, USA}
\email{xil0005@auburn.edu}
\thanks{Xiaolong Li was supported in part by the National Science
Foundation (DMS-2553660) and Auburn University start-up funding.}
\author{Elham Matinpour}
\address{Department of Mathematics and Statistics, Auburn University,
Auburn, AL 36849, USA}
\email{elm0135@auburn.edu}
\author{Peter Petersen}
\address{University of California, Los Angeles, 520 Portola Plaza,
Los Angeles, CA 90095, USA}
\email{petersen@math.ucla.edu}


\makeatother

\providecommand\theoremname{Theorem}
\theoremstyle{plain}

\begin{document}
\begin{abstract}
We prove that a closed connected Riemannian $n$-manifold, $n\geq4$,
satisfying $\sec>\frac{\Scal}{n(n+2)}$ at every point has vanishing
real cohomology in degrees $1,..., n-1$. It is consequently a rational
homology sphere if $n$ is odd or the manifold is orientable. Complex
projective space shows that the coefficient is optimal and that the
strict inequality is necessary in every even dimension. The required
estimate for the Weitzenböck curvature term is established by integrating
restricted Ricci operators over orthonormal frames and combining complementary
degrees using Hodge duality. This also allows us to replace sectional
curvature pinching by intermediate Ricci curvature bounds. For a nonflat
metric, the weak Yau inequality, $\sec\geq\frac{\Scal}{n(n+2)}$,
gives the same cohomology vanishing unless the manifold is isometric,
up to scaling, to a quotient of complex projective space with its
Fubini-Study metric. This also holds using the second intermediate
Ricci curvature in place of sectional curvature. More generally, a
nonzero Betti number in the range controlled by an intermediate Ricci
bound forces local symmetry. We also derive sharp Hodge eigenvalue
estimates and a scalar spectral vanishing criterion. 
\end{abstract}

\maketitle

\section{Introduction}

In 1990, Yau asked in his ``Open problems in geometry'' \cite[Problem~12, p.~4]{Yau}
(see also \cite[p.~369]{SY}): 
\begin{quote}
``The famous pinching problem says that on a compact simply connected
manifold if $K_{\min}>\frac{1}{4}K_{\max}>0$, then the manifold is
homeomorphic to a sphere. If we replace $K_{\max}$ by normalized
scalar curvature, can we deduce similar pinching theorems?'' 
\end{quote}
We write $\sec(\sigma)$ for the sectional curvature of a two-plane
$\sigma$. For a function $\lambda$ on $M$, the inequality $\sec\geq\lambda$
means $\sec(\sigma)\geq\lambda(x)$ for every $\sigma\subset T_{x}M$,
at every point $x$; strict and upper bounds have the analogous meaning.
When an extremum is needed, $\min\sec$ and $\max\sec$ are taken
over two-planes in a single tangent space. We use the same notation
for algebraic curvature tensors. The normalized scalar curvature is
$\frac{\Scal}{n(n-1)}$. Yau's pinching problem can be restated as
follows (see \cite[p.~525]{GX}).
\begin{problem}[Yau's pinching problem]
Let $(M^{n},g)$ be a closed simply connected Riemannian manifold
of dimension $n\geq4$. If 
\begin{equation}
\sec>\frac{\Scal}{n(n+2)}\label{eq:yau}
\end{equation}
holds, must $M$ be homeomorphic to $S^{n}$? 
\end{problem}

Complex projective space determines the coefficient in \eqref{eq:yau}.
In every even dimension $n=2m\geq4$, the Fubini-Study metric on $\CP^{m}$,
normalized by $1\leq\sec\leq4$, satisfies $\Scal=n(n+2)$. Since
$b_{2}(\CP^{m})=1$, strict inequality is necessary in every even
dimension.

We prove the following rational homology sphere theorem under Yau's
pinching condition in every dimension.

\begin{theorem}\label{thm:main} Let $(M^{n},g)$ be a closed connected
Riemannian manifold, $n\geq4$. If \eqref{eq:yau} holds, then 
\[
H^{j}(M;\R)=0\quad\text{for }1\leq j\leq n-1.
\]
If $n$ is odd or $M$ is orientable, then $M$ is a rational homology
sphere. \end{theorem}

Note that Theorem~\ref{thm:main} does not require simple connectivity.
However, it does not exclude torsion in integral homology, even when
$M$ is simply connected. No orientability hypothesis is needed for
the vanishing conclusion. The distinction in even dimensions concerns
only top-degree cohomology; see Section~\ref{sec:prelim}.

Earlier progress on Yau's problem includes both stronger pinching
theorems and results under additional geometric hypotheses. Gu and
Xu \cite[Theorem~1(ii)]{GX} proved that a closed $n$-manifold, $n\geq3$,
satisfying $\sec>\frac{\Scal}{n^{2}-n+6}$ is diffeomorphic to a spherical
space form. Costa and Ribeiro \cite[Corollary~1]{CR} settled the
simply connected 4-dimensional case through a weaker bound on biorthogonal
curvature. Xu and Gu \cite[Theorem~1.1 and Corollary~3.1]{XG} settled
the Einstein case and classified the simply connected Einstein metrics
with $\sec\geq\frac{\Scal}{n(n+2)}>0$ as round spheres or complex
projective spaces with the Fubini-Study metric. More recently, Li
\cite[Theorem~1.4]{Li} proved that $\sec>\frac{\Scal}{n^{2}-n+12}$
implies positive isotropic curvature and hence the topological sphere
conclusion in the simply connected case. Ge \cite[Theorem~1.5 and Corollary~1.6]{Ge}
strengthened the conclusion, without a fundamental-group assumption
and with the same coefficient $\frac{1}{n^{2}-n+12}$, to a diffeomorphic
spherical space form, by establishing positive complex sectional curvature
and applying Ricci flow. For $n>4$, these general sphere theorems
require a larger coefficient than Yau's.

With Yau's coefficient, Ge \cite[Theorem~1.2]{Ge} proved $b_{2}=b_{n-2}=0$.
In odd dimensions his sharp estimate on 2-forms works with the weaker
sufficient condition 
\[
\sec>\frac{n-3}{(n-1)(n^{2}-6)}\Scal.
\]
In the equality case, $\sec\geq\frac{\Scal}{n(n+2)}$, every harmonic
2-form is parallel. In the nonflat case, these forms vanish in odd
dimensions; in even dimensions, nonzero second cohomology forces the
manifold to be isometric, up to scaling, to complex projective space
with its Fubini-Study metric \cite[Theorem~1.3]{Ge}. Ge's result
already gives a rational homology sphere theorem in dimension five.
To obtain this conclusion in every dimension, we prove vanishing for
all degrees $1,...,n-1$.

Homology sphere theorems also arise from spectral conditions on curvature.
The classical work of Meyer \cite{Meyer} and Gallot-Meyer \cite{GM}
relates positivity and nonnegativity of the curvature operator to
the vanishing and parallelism of harmonic forms. Gallot \cite{Gallot}
subsequently obtained quantitative Betti-number bounds from lower
curvature-operator bounds and upper diameter bounds. Petersen and
Wink \cite[Corollary~to Theorem~A]{PW} proved that if the sum of
the lowest $\lceil\frac{n}{2}\rceil$ eigenvalues of the curvature
operator is positive, then all intermediate Betti numbers vanish.
A second approach uses the curvature operator of the second kind,
acting on trace-free symmetric two-tensors. Nienhaus, Petersen, and
Wink \cite[Theorem~A]{NPW} proved that a closed connected orientable
manifold with $\frac{n+2}{2}$-nonnegative curvature operator of the
second kind is flat or a rational homology sphere. Li \cite[Theorem~1.6]{LiCOSK}
enlarged this condition by allowing the average of the lowest $\frac{n+2}{2}$
eigenvalues to be negative, with a lower bound determined by their
overall mean. Li and Petersen \cite[Theorem~1.4]{LPCOSK} obtained
parallelism and vanishing under a corresponding partial-sum condition
after modifying this operator by Ricci and scalar curvature terms.
Here a partial sum includes the indicated fraction of the next eigenvalue.

The lower sectional and intermediate Ricci bounds in this paper and
these averaged spectral criteria do not, in general, imply each other.
Remark~\ref{rem:spectral-comparison} gives explicit seven-dimensional
algebraic curvature tensors showing that neither family of curvature
conditions implies the other.

We also obtain vanishing in a prescribed range of degrees under a
weaker curvature hypothesis. An algebraic curvature tensor $E$ has
\emph{nonnegative $p$th intermediate Ricci curvature} if 
\[
\sum^{p}_{i=1}E(v,e_{i},v,e_{i})\geq0
\]
for every orthonormal $(p+1)$-frame $v,e_{1},\ldots,e_{p}$; positivity
means that every such sum is strictly positive. We denote the minimum
of these sums by $\Ric_{p}(E)$; the precise definition and normalization
are given in Section~\ref{sec:prelim}. The corresponding lower bound
relative to scalar curvature gives the following result.

\begin{theorem}\label{thm:intermediate} Let $(M^{n},g)$ be a closed
connected Riemannian manifold, $n\geq4$, and let $1\leq p\leq\lfloor\frac{n}{2}\rfloor$
such that 
\begin{equation}
\frac{1}{p}\Ric_{p}(R)\geq\frac{\Scal}{n(n+2)}\label{eq:intermediate-pinching}
\end{equation}
at every point. It follows that every harmonic form of degree $p\leq j\leq n-p$
is parallel. If there is a point where \eqref{eq:intermediate-pinching}
is strict, then 
\[
b_{j}(M)=0\quad\text{for }p\leq j\leq n-p.
\]
\end{theorem}

Condition \eqref{eq:intermediate-pinching} says precisely that $R-\frac{\Scal}{n(n+2)}I$
has nonnegative $p$th intermediate Ricci curvature, where $I$ is
the tensor of constant sectional curvature one. For $p=1$ it is the
weak Yau inequality; for larger $p$ it bounds sums of sectional curvatures
and allows individual terms to fall below the same scalar threshold.
The weak inequality implies that harmonic forms are parallel, while
strictness at one point eliminates these parallel forms. The Fubini-Study
metric explains why the distinction is necessary.

The second intermediate Ricci bound gives vanishing in every intermediate
degree when it is strict at just one point.

\begin{corollary}\label{cor:second-ricci-strict} Let $(M^{n},g)$
be a closed connected Riemannian manifold, $n\geq4$. If 
\[
\frac{1}{2}\Ric_{2}(R)\geq\frac{\Scal}{n(n+2)}
\]
at every point, with strict inequality at some point, then $b_{j}(M)=0$
for $1\leq j\leq n-1$. \end{corollary}

No separate scalar curvature assumption is needed. Summing the two-term
inequalities and taking the trace implies 
\[
\Ric\geq\frac{n-1}{n(n+2)}\Scal\,g,\qquad\Scal\geq0,
\]
showing that Ricci tensor is nonnegative and positive at a point.
Theorem \ref{thm:algebraic} below offers a uniform proof of the
corollary in all degrees.

Intermediate Ricci curvature dates back at least to Bishop and Crittenden
\cite[pp.~253-257]{BC}, who considered these sums under the name
of $p$-mean curvature. Wu \cite[Theorem~5 and its Corollary]{Wu}
and Shen \cite[Theorem~1 and Corollary~1]{Shen} developed topological
consequences for complete open manifolds, including homology vanishing
under suitable hypotheses at infinity. Guijarro and Wilhelm established
focal-radius comparison and rigidity \cite{GW}, as well as connectivity
and sphere theorems under intermediate Ricci and conjugate-radius
bounds \cite{GWConnectivity}. Mouillé \cite{Mouille} proved symmetry
rank bounds and sphere theorems under positive second intermediate
Ricci curvature and torus symmetry assumptions. Complementing these
rigidity results, Reiser and Wraith \cite[Theorem~A]{RW} constructed
closed $n$-manifolds, $n\geq5$, with positive $k$th intermediate
Ricci curvature and arbitrarily large total Betti number for $\lfloor\frac{n}{2}\rfloor+2\leq k\leq n-1$.
These examples show that positive intermediate Ricci curvature alone
does not imply a general homology vanishing theorem.

We recall how these conditions compare with sectional and Ricci curvature.
Nonnegative first intermediate Ricci curvature is nonnegative sectional
curvature, and nonnegative $(n-1)$st intermediate Ricci curvature
is nonnegative Ricci curvature. Between these endpoints, the conditions
become weaker as the index increases: summing the $p$-term inequalities
over subsets gives every longer-sum inequality. They relate to partial
traces of the Jacobi operator, rather than sums of eigenvalues of
the Ricci tensor. Equivalently, nonnegative $p$th intermediate Ricci
curvature means that $\Ric(E|_{H})\geq0$ on every $(p+1)$-dimensional
subspace $H$. We use this description in the proof: Hodge duality
on $H$ identifies these restricted Ricci tensors with nonnegative
operators on $p$-forms.

We now describe the main step of the proof. Write $Q_{p}(R)$ for
the curvature term in the Bochner-Weitzenböck formula on $p$-forms,
normalized by $Q_{p}(I)=p(n-p)\Id$. The formula is 
\[
dd^{*}+d^{*}d=\nabla^{*}\nabla+Q_{p}(R).
\]
See \cite[Theorem~9.4.1]{Petersen} and \cite[formula~(1.154)]{Besse}.
For a harmonic $p$-form $\omega$ on a closed manifold, integration
by parts gives 
\begin{equation}
0=\int_{M}\bigl(|\nabla\omega|^{2}+\langle Q_{p}(R)\omega,\omega\rangle\bigr)\,dV.\label{eq:bochner}
\end{equation}
The integral uses the Riemannian volume density and requires no orientation.
The same conclusion follows from the identity 
\[
\frac{1}{2}\Delta|\omega|^{2}=-|\nabla\omega|^{2}-\langle Q_{p}(R)\omega,\omega\rangle,
\]
where $\Delta=d^{*}d$ on functions. If $Q_{p}(R)\geq0$, the maximum
principle makes $|\omega|$ constant, and the identity gives $\nabla\omega=0$.
If $Q_{p}(R)$ is also positive definite at one point, then $\omega=0$;
see \cite[1.137]{Besse}. Hodge theory gives $b_{p}=0$. The required
estimate is the following algebraic statement, whose hypothesis involves
only intermediate Ricci curvature.

\begin{theorem}\label{thm:algebraic} Let $n\geq4$ and $1\leq p\leq\lfloor\frac{n}{2}\rfloor$.
If an algebraic curvature tensor $E$ on $\R^{n}$ has nonnegative
$p$th intermediate Ricci curvature, then 
\begin{equation}
Q_{p}(E)\geq-\frac{p(n-p)}{3n}\Scal(E)\Id.\label{eq:algebraic}
\end{equation}
The inequality is strict if its $p$th intermediate Ricci curvature
is positive. Moreover, for every algebraic curvature tensor $R$ and
every $1\leq p\leq n-1$, 
\begin{equation}
Q_{p}(R)\geq\frac{p(n-p)}{3n}\left(\frac{n(n+2)}{2}\Ric_{2}(R)-\Scal(R)\right)\Id.\label{eq:main-estimate}
\end{equation}
\end{theorem}

To prove Theorem~\ref{thm:algebraic}, we average the operators $Q_{p}(E|_{H})$
over all coordinate $(p+1)$-planes in all orthonormal frames, extending
each operator by zero outside $\Lambda^{p}H^{*}$. Hodge duality within
$H$ identifies it with $\Ric(E|_{H})$, so the average is nonnegative
under the hypothesis of Theorem~\ref{thm:algebraic}. We compute
the average as a linear combination of $Q_{p}(E)$, the derivation
action of $\Ric(E)$, and scalar curvature. Since $n-p\geq p$, the
same positivity holds in the complementary degree. Hodge duality on
the ambient space then allows the two identities to be combined so
that their Ricci terms cancel. The result is \eqref{eq:algebraic},
with a coefficient that produces the same curvature pinching in every
degree.

Applying the same identity to $R-\frac{1}{2}\Ric_{2}(R)I$ gives \eqref{eq:main-estimate}
in every degree, including degrees 1 and $n-1$. This proves Corollary~\ref{cor:second-ricci-strict}
and Theorem~\ref{thm:main}. Applying it to $R-\frac{\Scal}{n(n+2)}I$
proves Theorem~\ref{thm:intermediate}. Thus the same complementary-degree
cancellation explains the common curvature pinching and all the vanishing
conclusions.

Invariant integration has a long history in geometry. Berger expressed
the scalar curvature of a Kähler metric as an integral of its holomorphic
sectional curvature \cite{Berger1966}; see also \cite[Corollary~3]{Magnusson}.
Gray used integration on the unit tangent sphere bundle to study compact
Kähler manifolds with constant scalar curvature and nonnegative sectional
curvature \cite{Gray}; a modern account is \cite[Section~2]{MSW}.
Bettiol and Mendes \cite[Proposition~3.1]{BM} express the Weitzenböck
curvature term on trace-free symmetric tensors as a sphere integral
involving sectional curvature. We apply invariant integration to the
Ricci operators of subspaces and combine the identities in complementary
degrees.

For weak curvature pinching, parallel harmonic forms lead to rigidity.
The curvature estimate makes harmonic forms parallel in the controlled
degrees, while the pinching condition forces irreducible restricted
holonomy and excludes Ricci flatness in the nonflat case. For the
second intermediate Ricci bound, harmonic one-forms also vanish, and
a restricted-root argument excludes symmetric spaces of higher rank.
The remaining cases are settled by the holonomy classification and
the Kähler and quaternionic curvature identities proved in Section~\ref{sec:odd}.
We obtain the following extension of Theorem~\ref{thm:main}.

\begin{theorem}\label{thm:odd-weak} Let $(M^{n},g)$ be a closed
connected Riemannian manifold, $n\geq4$. Suppose that $g$ is not
flat and satisfies 
\begin{equation}
\frac{1}{2}\Ric_{2}(R)\geq\frac{\Scal}{n(n+2)}\label{eq:second-ricci-pinching}
\end{equation}
at every point. It follows that either 
\[
H^{j}(M;\R)=0\quad\text{for }1\leq j\leq n-1,
\]
or $n$ is even and $(M,g)$ is isometric, up to scaling, to a quotient
of $\CP^{\frac{n}{2}}$ with its Fubini-Study metric. In particular,
$M$ is a rational homology sphere when $n$ is odd. \end{theorem}

If $M$ is orientable, the first alternative gives a rational homology
sphere, and the complex-projective quotient in the second is eliminated
by Synge's theorem. For general $p$, the same holonomy argument shows
that nonzero cohomology in the controlled range forces local symmetry.

\begin{theorem}\label{thm:intermediate-rigidity} Let $(M^{n},g)$
be a closed connected Riemannian manifold, $n\geq4$, and let $1\leq p\leq\lfloor\frac{n}{2}\rfloor$.
Suppose that $g$ is not flat and satisfies \eqref{eq:intermediate-pinching}
at every point. If $g$ is not locally symmetric, then 
\[
b_{j}(M)=0\quad\text{for }p\leq j\leq n-p.
\]
Consequently, if $b_{j}(M)\ne0$ for some $p\leq j\leq n-p$, then
$g$ is locally symmetric and its universal Riemannian cover is a
compact irreducible symmetric space. \end{theorem}

In the companion paper \cite{LMPUpper}, we apply the same integral
identity to upper bounds for sectional curvature relative to scalar
curvature. That paper proves homology sphere theorems and determines
the sharp coefficients and optimal pinching in dimensions four and
five.

The paper is organized as follows. In Section~\ref{sec:prelim},
we fix conventions and recall the sphere and frame integration formulas.
In Section~\ref{sec:integral}, we prove the integral identity in
Proposition~\ref{prop:average}. We apply it in Section~\ref{sec:proof}
to prove the curvature estimate and the vanishing results, and to
discuss sharpness in even dimensions. In Section~\ref{sec:odd},
we prove Theorems~\ref{thm:odd-weak} and \ref{thm:intermediate-rigidity}
and record the improvement in dimension five. Finally, in Section~\ref{sec:spectral-pinching},
we prove Hodge eigenvalue estimates and a vanishing theorem under
a scalar spectral condition.

\medskip{}
\noindent\textbf{AI disclosure.} The authors used OpenAI's ChatGPT
as an assistive tool in the development of the paper. The authors
take full responsibility for all mathematical arguments.

\section{Preliminaries}

\label{sec:prelim}

Throughout the paper, manifolds are connected. We recall two facts
that allow us to work without an orientability assumption. If $\pi:\widehat{M}\to M$
is a finite covering with $r$ sheets, summing pullbacks over local
inverse images of $\pi$ defines the transfer $\pi_{*}:\Omega^{p}\left(\widehat{M}\right)\to\Omega^{p}\left(M\right)$.
It commutes with exterior differentiation and satisfies 
\[
\pi_{*}\pi^{*}=r\Id.
\]
Consequently, pullback is injective on de Rham cohomology: 
\[
\pi^{*}:H^{j}(M;\R)\longrightarrow H^{j}(\widehat{M};\R).
\]
We may therefore prove cohomology vanishing on the orientation double
cover when necessary. For a Riemannian covering, pullback also preserves
harmonic and parallel forms.

If all real cohomology groups in degrees $1,...,n-1$ of a closed
manifold vanish, then orientability gives the rational homology of
a sphere. The Euler characteristic satisfies $\chi\left(\hat{M}\right)=r\chi\left(M\right)$
for an $r$-fold covering space. Poincaré duality on the orientation
cover therefore shows that $\chi(M)=0$ in odd dimensions. In that
case, vanishing of the intermediate Betti numbers shows that 
\[
0=\chi(M)=1-b_{n}(M),
\]
proving that $b_{n}(M)=1$. Thus our vanishing conclusions give rational
homology spheres whenever the dimension is odd or the manifold is
orientable.

Let $(V,g)$ be an $n$-dimensional Euclidean space. The inner product
$g$ identifies $V$ with its dual $V^{*}$ by the isomorphism $v\mapsto g(v,\cdot)$.
We use the induced identifications of their exterior powers and identify
symmetric two-tensors with self-adjoint endomorphisms through $g$.
For an orthonormal basis $\{e_{1},\ldots,e_{n}\}$, the corresponding
dual basis is $\{e^{1},\ldots,e^{n}\}$, where $e^{i}=g(e_{i},\cdot)$.

\subsection{Algebraic curvature tensors}

An \emph{algebraic curvature tensor} on $V$ is a multilinear map
$R\colon V^{4}\to\R$ satisfying the symmetries 
\[
R(X,Y,Z,W)=-R(Y,X,Z,W)=-R(X,Y,W,Z)=R(Z,W,X,Y)
\]
and the first Bianchi identity 
\[
R(X,Y,Z,W)+R(X,Z,W,Y)+R(X,W,Y,Z)=0
\]
for all $X,Y,Z,W\in V$. Write $R_{ijkl}=R(e_{i},e_{j},e_{k},e_{l})$.
For orthonormal $X,Y$, the sectional curvature of the two-plane spanned
by $X$ and $Y$ is 
\[
\sec\bigl(\operatorname{span}\{X,Y\}\bigr)=R(X,Y,X,Y).
\]
The Ricci and scalar curvatures are the contractions 
\[
\Ric(R)(X,Y)=\sum^{n}_{a=1}R(X,e_{a},Y,e_{a})
\]
and 
\[
\Scal(R)=\sum^{n}_{a=1}\Ric(R)(e_{a},e_{a})=2\sum_{a<b}R_{abab}.
\]
These expressions are independent of the orthonormal basis.

For $1\leq p\leq n-1$, define the \emph{minimum $p$th intermediate
Ricci curvature} by 
\[
\Ric_{p}(R)=\min_{\substack{|v|=1,\ H\subset v^{\perp}\\
\dim H=p
}
}\sum^{p}_{i=1}R(v,e_{i},v,e_{i}),
\]
where $e_{1},\ldots,e_{p}$ is any orthonormal basis of $H$. The
sum is independent of this basis: it is the partial trace on $H$
of the Jacobi operator defined by $\langle J_{v}x,y\rangle=R(x,v,y,v)$.
We use the unnormalized convention. Thus $\Ric_{p}(R)\geq c$ means
\[
\sum^{p}_{i=1}R(v,e_{i},v,e_{i})\geq c
\]
for every orthonormal frame $v,e_{1},\ldots,e_{p}$, and similarly
for strict inequalities. The extreme cases are 
\[
\Ric_{1}(R)=\min\sec(R)
\]
and 
\[
\Ric_{n-1}(R)=\min_{|v|=1}\Ric(R)(v,v).
\]
Equivalently, $\Ric_{p}(R)$ is the minimum over unit $v$ of the
sum of the $p$ smallest eigenvalues of $J_{v}|_{v^{\perp}}$.

The tensor $I$ of constant sectional curvature 1 is 
\[
I(X,Y,Z,W)=\langle X,Z\rangle\langle Y,W\rangle-\langle X,W\rangle\langle Y,Z\rangle.
\]
In particular, $\Ric_{p}(I)=p$, and for every $t\in\R$, 
\[
\Ric_{p}(R-tI)=\Ric_{p}(R)-pt.
\]
For a symmetric two-tensor $L$, its Kulkarni-Nomizu product with
$g$ is 
\[
(L\owedge g)_{ijkl}=L_{ik}\delta_{jl}+L_{jl}\delta_{ik}-L_{il}\delta_{jk}-L_{jk}\delta_{il}.
\]
Thus $g\owedge g=2I$, and direct contraction gives 
\[
\Ric(L\owedge g)=(n-2)L+(\tr L)g.
\]
For $n\geq3$, every algebraic curvature tensor has the decomposition
\[
R=\frac{\Scal(R)}{n(n-1)}I+\frac{1}{n-2}\left(\Ric(R)-\frac{\Scal(R)}{n}g\right)\owedge g+W.
\]
The three terms are the scalar curvature component, the trace-free
Ricci component, and the Weyl tensor $W$, respectively.

We also regard $R$ as the self-adjoint operator on $\Lambda^{2}V^{*}$
determined by 
\[
\langle R(e^{i}\wedge e^{j}),e^{k}\wedge e^{l}\rangle=R_{ijkl}.
\]
The inner product and norm on curvature tensors are those of these
operators. Thus 
\[
|R|^{2}=\sum_{i<j,\,k<l}R^{2}_{ijkl}=\frac{1}{4}\sum_{i,j,k,l}R^{2}_{ijkl}.
\]
With this convention, the Ricci contraction is adjoint to the Kulkarni-Nomizu
product with $g$, namely 
\begin{equation}
\langle R,L\owedge g\rangle=\langle\Ric(R),L\rangle.\label{eq:ricci-adjoint}
\end{equation}
Indeed, the curvature symmetries give 
\begin{align*}
\langle R,L\owedge g\rangle & =\frac{1}{4}\sum_{i,j,k,l}R_{ijkl}(L_{ik}\delta_{jl}+L_{jl}\delta_{ik}-L_{il}\delta_{jk}-L_{jk}\delta_{il})\\
 & =\sum_{i,k}\left(\sum_{j}R_{ijkj}\right)L_{ik}=\langle\Ric(R),L\rangle.
\end{align*}
In particular, the scalar, trace-free Ricci, and Weyl components are
pairwise orthogonal. The identity endomorphism on any exterior power
is denoted by $\Id$.

\subsection{The derivation action}

For a linear endomorphism $L$ of $V^{*}$, its \emph{derivation action}
is the unique extension to the exterior algebra $\Lambda^{*}V^{*}$
that preserves degree, vanishes on scalars, and satisfies the Leibniz
rule. Writing $L^{[p]}$ for its restriction to $\Lambda^{p}V^{*}$,
this means $L^{[0]}=0$, $L^{[1]}=L$, and 
\[
L^{[r+s]}(\alpha\wedge\beta)=\bigl(L^{[r]}\alpha\bigr)\wedge\beta+\alpha\wedge\bigl(L^{[s]}\beta\bigr)
\]
for $\alpha\in\Lambda^{r}V^{*}$ and $\beta\in\Lambda^{s}V^{*}$.
Equivalently, on decomposable $p$-forms, 
\[
L^{[p]}(\alpha_{1}\wedge\cdots\wedge\alpha_{p})=\sum^{p}_{r=1}\alpha_{1}\wedge\cdots\wedge L\alpha_{r}\wedge\cdots\wedge\alpha_{p}.
\]
Thus $L^{[p]}$ applies $L$ to one factor at a time. If $L$ is self-adjoint,
so is $L^{[p]}$.

For $i<j$, identify the bivector $e_{i}\wedge e_{j}$ with the skew-adjoint
endomorphism of $V$ given by 
\[
(e_{i}\wedge e_{j})v=g(e_{i},v)e_{j}-g(e_{j},v)e_{i}.
\]
Transport this endomorphism to $V^{*}$ through the metric identification
and extend it to exterior forms by its derivation action, using the
same notation. Thus $(e_{i}\wedge e_{j})e^{i}=e^{j}$, $(e_{i}\wedge e_{j})e^{j}=-e^{i}$,
and the other basis covectors are annihilated. On any exterior form
$\omega$, the action is 
\[
(e_{i}\wedge e_{j})\omega=e^{j}\wedge\iota_{e_{i}}\omega-e^{i}\wedge\iota_{e_{j}}\omega,
\]
where $\iota_{v}$ denotes contraction by $v$: 
\[
(\iota_{v}\omega)(v_{2},\ldots,v_{p})=\omega(v,v_{2},\ldots,v_{p}),\quad\text{for }\omega\in\Lambda^{p}V^{*},\ p\geq1.
\]
It is zero on scalars. Wedge multiplication by $e^{i}$ is adjoint
to $\iota_{e_{i}}$, so each infinitesimal rotation is skew-adjoint
on every exterior power.

The Weitzenböck curvature term on $\Lambda^{p}V^{*}$ is 
\begin{equation}
Q_{p}(R)=-\sum_{i<j,\,k<l}R_{ijkl}(e_{i}\wedge e_{j})\circ(e_{k}\wedge e_{l}).\label{eq:q-definition}
\end{equation}
Here $\circ$ denotes composition of the derivation actions, each
of which preserves exterior degree. This definition is independent
of the orthonormal basis and is linear in $R$. The curvature symmetry
$R_{ijkl}=R_{klij}$ also gives 
\[
Q_{p}(R)^{*}=-\sum_{i<j,\,k<l}R_{ijkl}(e_{k}\wedge e_{l})\circ(e_{i}\wedge e_{j})=Q_{p}(R).
\]
See \cite[Section~1.I]{Besse} for the general construction and \cite[formula~(2.2)]{BM}
for this normalization.

\begin{lemma}\label{lem:q-identities} Let $R$ be an algebraic curvature
tensor on $V$, let $L$ be a symmetric 2-tensor, regarded as a self-adjoint
endomorphism, and let $1\leq p\leq n-1$. It follows that 
\begin{equation}
Q_{1}(R)=\Ric(R).\label{eq:q-ricci}
\end{equation}
\begin{equation}
Q_{p}(I)=p(n-p)\Id.\label{eq:q-constant}
\end{equation}
\begin{equation}
*^{-1}Q_{n-p}(R)*=Q_{p}(R).\label{eq:q-hodge}
\end{equation}
\begin{equation}
Q_{p}(L\owedge g)=(n-2p)L^{[p]}+p(\tr L)\Id.\label{eq:kulkarni-q}
\end{equation}
Either orientation can be used in \eqref{eq:q-hodge}. \end{lemma}
\begin{proof}
Applying \eqref{eq:q-definition} to a basis covector gives 
\[
Q_{1}(R)e^{j}=\sum_{i}\left(\sum_{a}R_{iaja}\right)e^{i},
\]
which proves \eqref{eq:q-ricci}.

To prove \eqref{eq:kulkarni-q}, diagonalize $L$ with eigenvalues
$\lambda_{1},\ldots,\lambda_{n}$. The curvature operator $L\owedge g$
is then diagonal in the bivector basis, with eigenvalues $\lambda_{i}+\lambda_{j}$.
For a $p$-element index set $A$, write $e^{A}$ for the wedge of
the corresponding basis covectors in increasing order. In \eqref{eq:q-definition},
$-\left(e_{i}\wedge e_{j}\right)\circ\left(e_{i}\wedge e_{j}\right)$
acts as the identity on $e^{A}$ when exactly one of $i,j$ belongs
to $A$, and as zero otherwise. Thus $Q_{p}(L\owedge g)$ acts on
$e^{A}$ with eigenvalue 
\begin{align*}
\sum_{i\in A,\,a\notin A}(\lambda_{i}+\lambda_{a}) & =(n-p)\sum_{i\in A}\lambda_{i}+p\sum_{a\notin A}\lambda_{a}\\
 & =(n-2p)\sum_{i\in A}\lambda_{i}+p\,\tr L.
\end{align*}
This is the eigenvalue of the right-hand side of \eqref{eq:kulkarni-q}.
Taking $L=\frac{1}{2}g$, so that $L\owedge g=I$, also proves \eqref{eq:q-constant}.

Finally, the Hodge star commutes with every orientation-preserving
orthogonal transformation. Differentiating along a rotation in the
$e_{i},e_{j}$-plane gives 
\[
*((e_{i}\wedge e_{j})\omega)=(e_{i}\wedge e_{j})(*\omega).
\]
It therefore commutes with the compositions in \eqref{eq:q-definition},
which proves \eqref{eq:q-hodge}. 
\end{proof}

For self-adjoint $L$, diagonalization also implies 
\begin{equation}
*^{-1}L^{[n-p]}*=(\tr L)\Id-L^{[p]}.\label{eq:hodge-derivation}
\end{equation}
Indeed, in an eigenbasis of $L$, the form $*e^{A}$ has the complementary
index set, and hence 
\[
L^{[n-p]}(*e^{A})=\left(\sum_{a\notin A}\lambda_{a}\right)*e^{A}=\left(\tr L-\sum_{i\in A}\lambda_{i}\right)*e^{A}.
\]

\subsection{Spherical integration}

For a Euclidean space $U$, write $S(U)$ for its unit sphere. On
each such sphere, $d\sigma$ denotes its Riemannian measure divided
by its total area. On the space of orthonormal frames of $V$, identified
with $\mathrm{O}(n)$ after fixing one frame, $d\mathcal{B}$ denotes
the invariant volume measure normalized by $\int_{\mathrm{O}(n)}d\mathcal{B}=1$.
All integrals of tensors or endomorphisms are taken componentwise
in the fixed ambient vector space. Orthogonal invariance and these
normalizations give 
\begin{equation}
\int_{\mathrm{O}(n)}\sum^{n}_{i=1}f(e_{i})\,d\mathcal{B}=n\int_{S(V)}f(u)\,d\sigma(u),\quad\mathcal{B}=(e_{1},\ldots,e_{n}),\label{eq:frame-sphere}
\end{equation}
for every continuous function $f$, including vector-valued functions.
Indeed, for each $i$, integration of $f(e_{i})$ over frames induces
the rotation-invariant integral of $f$ over the sphere with total
mass one.

\begin{lemma}[Spherical integration]\label{lem:sphere-integrals}
If $U$ has dimension $d\geq2$ and $u_{1},\ldots,u_{d}$ are orthonormal
coordinates, then 
\begin{equation}
\int_{S(U)}u_{i}u_{j}\,d\sigma(u)=\frac{\delta_{ij}}{d}\label{eq:sphere-quadratic}
\end{equation}
and 
\begin{equation}
\int_{S(U)}u_{i}u_{j}u_{k}u_{l}\,d\sigma(u)=\frac{1}{d(d+2)}\bigl(\delta_{ij}\delta_{kl}+\delta_{ik}\delta_{jl}+\delta_{il}\delta_{jk}\bigr).\label{eq:sphere-quartic}
\end{equation}
\end{lemma}
\begin{proof}
Reflection in one coordinate shows that an integral vanishes if that
coordinate occurs to an odd power. Permutation invariance and $\sum_{i}u^{2}_{i}=1$
then give 
\[
\int_{S(U)}u^{2}_{i}\,d\sigma(u)=\frac{1}{d}\int_{S(U)}\sum^{d}_{a=1}u^{2}_{a}\,d\sigma(u)=\frac{1}{d},
\]
which proves \eqref{eq:sphere-quadratic}.

For the fourth-order identity, rotate the first two coordinates through
$\frac{\pi}{4}$. Expanding and discarding the odd terms gives 
\begin{align*}
\int_{S(U)}u^{4}_{1}\,d\sigma(u) & =\frac{1}{4}\int_{S(U)}(u_{1}+u_{2})^{4}\,d\sigma(u)\\
 & =\frac{1}{2}\int_{S(U)}u^{4}_{1}\,d\sigma(u)+\frac{3}{2}\int_{S(U)}u^{2}_{1}u^{2}_{2}\,d\sigma(u).
\end{align*}
Thus $\int u^{4}_{1}\,d\sigma=3\int u^{2}_{1}u^{2}_{2}\,d\sigma$.
Integrating $\bigl(\sum_{i}u^{2}_{i}\bigr)^{2}=1$ now yields 
\begin{align*}
1 & =d\int_{S(U)}u^{4}_{1}\,d\sigma(u)+d(d-1)\int_{S(U)}u^{2}_{1}u^{2}_{2}\,d\sigma(u)\\
 & =d(d+2)\int_{S(U)}u^{2}_{1}u^{2}_{2}\,d\sigma(u).
\end{align*}
These values, together with reflection and permutation invariance,
give all components of \eqref{eq:sphere-quartic}. 
\end{proof}

\section{The integral identity}
\label{sec:integral}

Let $R$ be an algebraic curvature tensor on an $n$-dimensional Euclidean
space $V$, and let $1\leq p\leq n-1$. For each $(p+1)$-dimensional
subspace $H$, extend $Q_{p}(R|_{H})$ by zero from $\Lambda^{p}H^{*}$
to $\Lambda^{p}V^{*}$. For an orthonormal frame $\mathcal{B}=(e_{1},\ldots,e_{n})$,
set 
\[
B^{(p)}_{\mathcal{B}}(R)=\sum_{\substack{H=\operatorname{span}\left\{ e_{i}\mid i\in A\right\} \\
|A|=p+1}}Q_{p}(R|_{H})
\]
and 
\[
\overline{B}_{p}(R)=\int_{\mathrm{O}(n)}B^{(p)}_{\mathcal{B}}(R)\,d\mathcal{B}.
\]

\begin{proposition}\label{prop:average} As self-adjoint operators
on $\Lambda^{p}V^{*}$, we have 
\begin{equation}
(n+2)\overline{B}_{p}(R)=3Q_{p}(R)-p\Ric(R)^{[p]}+p\Scal(R)\Id.\label{eq:average}
\end{equation}
\end{proposition}

Before proving the identity, we record the cancellation that explains
its role in the pinching theorems. Applying \eqref{eq:average} in
degree $n-p$ and conjugating by the Hodge star, equations \eqref{eq:q-hodge}
and \eqref{eq:hodge-derivation} give 
\[
(n+2)*^{-1}\overline{B}_{n-p}(R)*=3Q_{p}(R)+(n-p)\Ric(R)^{[p]}.
\]
Multiply \eqref{eq:average} by $n-p$ and this identity by $p$.
The Ricci terms cancel, leaving 
\begin{equation}
3nQ_{p}(R)+p(n-p)\Scal(R)\Id=(n+2)\bigl((n-p)\overline{B}_{p}(R)+p*^{-1}\overline{B}_{n-p}(R)*\bigr).\label{eq:positive-representation}
\end{equation}
Thus nonnegativity of the two averages gives \eqref{eq:algebraic}.
For the shifted tensor $R-tI$, the left-hand side becomes 
\[
3nQ_{p}(R)+p(n-p)\bigl(\Scal(R)-n(n+2)t\bigr)\Id.
\]
The choice $t=\frac{\Scal(R)}{n(n+2)}$ cancels the scalar term in
every degree. This is the source of the common pinching condition.

We now prove Proposition~\ref{prop:average}, first in a fixed frame
and then by integration. Let $P_{a}$ be the orthogonal projection
of $\Lambda^{2}V^{*}$ onto 
\[
e^{a}\wedge(e^{a})^{\perp}=\operatorname{span}\{e^{a}\wedge e^{j}:j\ne a\},
\]
and define 
\[
T_{\mathcal{B}}R=\frac{1}{2}\sum^{n}_{a=1}P_{a}\circ R\circ P_{a}.
\]
Since $P_{a}(e^{i}\wedge e^{j})=(\delta_{ia}+\delta_{ja})e^{i}\wedge e^{j}$,
its components are 
\[
(T_{\mathcal{B}}R)_{ijkl}=\frac{1}{2}\sum_{a}(\delta_{ia}+\delta_{ja})(\delta_{ka}+\delta_{la})R_{ijkl}.
\]
Thus $T_{\mathcal{B}}$ multiplies curvature components with two,
three, or four distinct indices by $1$, $\frac{1}{2}$, or $0$,
respectively. These multiples preserve the curvature symmetries and
the Bianchi identity, so $T_{\mathcal{B}}R$ is an algebraic curvature
tensor. The same component formula shows that $T_{\mathcal{B}}$ is
self-adjoint for the Hilbert-Schmidt inner product on algebraic curvature
tensors.

\begin{lemma}\label{lem:coordinate-identity} For every algebraic
curvature tensor $R$, every orthonormal frame $\mathcal{B}$, and
$1\leq p\leq n-1$, 
\begin{equation}
B^{(p)}_{\mathcal{B}}(R)=Q_{p}(T_{\mathcal{B}}R)-\frac{1}{2}\bigl(\Ric(R)-\Ric(R)_{\mathrm{d}}\bigr)^{[p]}.\label{eq:direct-coordinate}
\end{equation}
Here $\Ric(R)_{\mathrm{d}}=\sum_{i}\Ric(R)_{ii}e_{i}\otimes e_{i}$
is the diagonal part of $\Ric(R)$ in $\mathcal{B}$. \end{lemma}
\begin{proof}
Write $e^{A}$ for the wedge of the basis covectors indexed by a $p$-element
set $A$, in increasing order. We compare matrix entries according
to how many indices the two forms have in common. Empty wedges equal
$1$ and empty sums equal $0$. The extension by zero ensures that
a coordinate plane $H$ can contribute to $\langle Q_{p}(R|_{H})e^{C},e^{A}\rangle$
only if it contains all directions indexed by $A\cup C$.

As in the proof of Lemma~\ref{lem:q-identities}, a rotation $e_{i}\wedge e_{j}$
sends $e^{A}$ to a signed unit basis form when exactly one of $i,j$
belongs to $A$, and annihilates it otherwise. The nonzero images
for distinct pairs are orthogonal. Since the rotations are skew-adjoint,
\eqref{eq:q-definition} gives 
\begin{equation}
\langle Q_{p}(R)e^{A},e^{A}\rangle=\sum_{i\in A,\,a\notin A}R_{iaia}.\label{eq:q-diagonal}
\end{equation}
For the diagonal entry of $B^{(p)}_{\mathcal{B}}(R)$, only the planes
$H=\operatorname{span}\{e_{a},e_{i}:i\in A\}$ with $a\notin A$ contribute.
Orient $H$ by $e^{a}\wedge e^{A}$, so $*_{H}e^{a}=e^{A}$. Hodge
duality on $H$ gives 
\[
\langle Q_{p}(R|_{H})e^{A},e^{A}\rangle=\Ric(R|_{H})(e_{a},e_{a})=\sum_{i\in A}R_{iaia}.
\]
Summing over $a\notin A$ gives the right-hand side of \eqref{eq:q-diagonal}.
The tensor $T_{\mathcal{B}}R$ has the same sectional components as
$R$, and the Ricci correction in \eqref{eq:direct-coordinate} has
zero diagonal. This proves the identity on diagonal entries.

Suppose next that the forms differ in one index. Write them as $e^{i}\wedge e^{S}$
and $e^{j}\wedge e^{S}$, where $|S|=p-1$ and $i,j\notin S$ are
distinct. The displayed order fixes their signs. Let $O$ be the complement
of $S\cup\{i,j\}$. Only $H=\operatorname{span}\{e_{i},e_{j},e_{s}:s\in S\}$
can contribute. With orientation $e^{i}\wedge e^{j}\wedge e^{S}$,
its Hodge star satisfies 
\[
*_{H}e^{i}=e^{j}\wedge e^{S},\qquad*_{H}e^{j}=-e^{i}\wedge e^{S}.
\]
Hence 
\begin{equation}
\langle B^{(p)}_{\mathcal{B}}(R)(e^{j}\wedge e^{S}),e^{i}\wedge e^{S}\rangle=-\Ric(R|_{H})_{ij}=-\sum_{s\in S}R_{isjs}.\label{eq:average-offdiagonal}
\end{equation}
Polarizing \eqref{eq:q-diagonal} in $u_{\pm}=\frac{1}{\sqrt{2}}(e_{i}\pm e_{j})$
gives 
\[
\langle Q_{p}(R)(e^{j}\wedge e^{S}),e^{i}\wedge e^{S}\rangle=\sum_{a\in O}R_{iaja}-\sum_{s\in S}R_{isjs}.
\]
Indeed, the terms involving $u_{\pm}$ in the form change by $2R_{iaja}$
for $a\in O$, whereas those involving $u_{\mp}$ in its orthogonal
complement change by $-2R_{isjs}$ for $s\in S$; all other terms
cancel. Every component in this formula has three distinct indices.
Therefore, we have 
\[
\langle Q_{p}(T_{\mathcal{B}}R)(e^{j}\wedge e^{S}),e^{i}\wedge e^{S}\rangle  =\frac{1}{2}\left(\sum_{a\in O}R_{iaja}-\sum_{s\in S}R_{isjs}\right)
\]
and 
\[
\langle(\Ric(R)-\Ric(R)_{\mathrm{d}})^{[p]}(e^{j}\wedge e^{S}),e^{i}\wedge e^{S}\rangle  =\Ric(R)_{ij}=\sum_{a\in O}R_{iaja}+\sum_{s\in S}R_{isjs}.
\]
Subtracting half of the second entry from the first gives $-\sum_{s\in S}R_{isjs}$,
as required by \eqref{eq:average-offdiagonal}.

Finally, if the forms differ in at least two indices, then $|A\cup C|\geq p+2$.
No coordinate $(p+1)$-plane contains their union, so the entry of
$B^{(p)}_{\mathcal{B}}(R)$ is zero. The Ricci derivation also has
zero entry, since it replaces at most one index. For a composition
of two rotations to connect the forms, their defining pairs must contain
all indices in $(A\setminus C)\cup(C\setminus A)$. This requires
at least four distinct indices, for which the corresponding component
of $T_{\mathcal{B}}R$ is zero. Thus the entry of $Q_{p}(T_{\mathcal{B}}R)$
vanishes as well, completing the proof. 
\end{proof}

\begin{proof}[Proof of Proposition~\ref{prop:average}]
By Lemma~\ref{lem:coordinate-identity}, it remains to average $T_{\mathcal{B}}R$
and $\Ric(R)_{\mathrm{d}}$ over frames. Recall the Jacobi operator
$\langle J_{u}v,w\rangle=R(v,u,w,u)$. It is self-adjoint and annihilates
$u$. The isometry $v\mapsto e_{a}\wedge v$ from $e^{\perp}_{a}$
onto $\operatorname{im}P_{a}$ identifies $J_{e_{a}}|_{e^{\perp}_{a}}$
with $P_{a}RP_{a}|_{\operatorname{im}P_{a}}$. Consequently, 
\[
\langle P_{a}RP_{a},R\rangle=\tr\left(P^{2}_{a}RP^{2}_{a}R\right)=\tr\bigl((P_{a}RP_{a})^{2}\bigr)=\tr(J^{2}_{e_{a}}),
\]
and hence 
\begin{equation}
\langle T_{\mathcal{B}}R,R\rangle=\frac{1}{2}\sum^{n}_{a=1}\tr(J^{2}_{e_{a}}).\label{eq:direct-jacobi-frame}
\end{equation}

In a fixed orthonormal basis, $\langle J_{u}e_{i},e_{j}\rangle=\sum_{a,b}R_{iajb}u_{a}u_{b}$.
The quartic spherical integral \eqref{eq:sphere-quartic} therefore
gives 
\[
\int_{S(V)}\tr(J^{2}_{u})\,d\sigma(u)=\frac{1}{n(n+2)}\sum_{i,j,a,b}\bigl(R_{iaja}R_{ibjb}+R^{2}_{iajb}+R_{iajb}R_{ibja}\bigr).
\]
The first two contractions are $|\Ric(R)|^{2}$ and $4|R|^{2}$. For
the third, the Bianchi identity $R_{iajb}-R_{ibja}=R_{ijab}$ yields
\[
4|R|^{2}=\sum_{i,j,a,b}(R_{iajb}-R_{ibja})^{2}=8|R|^{2}-2\sum_{i,j,a,b}R_{iajb}R_{ibja}.
\]
Thus 
\begin{equation}
\int_{S(V)}\tr(J^{2}_{u})\,d\sigma(u)=\frac{1}{n(n+2)}\bigl(|\Ric(R)|^{2}+6|R|^{2}\bigr).\label{eq:direct-jacobi-sphere}
\end{equation}
Combining \eqref{eq:frame-sphere}, \eqref{eq:direct-jacobi-frame},
and \eqref{eq:direct-jacobi-sphere}, we obtain 
\[
\int_{\mathrm{O}(n)}\langle T_{\mathcal{B}}R,R\rangle\,d\mathcal{B}=\frac{1}{n+2}\left(3|R|^{2}+\frac{1}{2}|\Ric(R)|^{2}\right).
\]
Since $T_{\mathcal{B}}$ is self-adjoint, polarization and the adjoint
identity \eqref{eq:ricci-adjoint} now give 
\begin{equation}
\int_{\mathrm{O}(n)}T_{\mathcal{B}}R\,d\mathcal{B}=\frac{1}{n+2}\left(3R+\frac{1}{2}\Ric(R)\owedge g\right).\label{eq:direct-curvature-integral}
\end{equation}
Applying the same quartic spherical formula to the diagonal Ricci
tensor gives 
\begin{align}
\int_{\mathrm{O}(n)}\Ric(R)_{\mathrm{d}}\,d\mathcal{B} & =n\int_{S(V)}\Ric(R)(u,u)\,u\otimes u\,d\sigma(u)\nonumber \\
 & =\frac{1}{n+2}\bigl(2\Ric(R)+\Scal(R)g\bigr).\label{eq:direct-ricci-integral}
\end{align}
Integrating \eqref{eq:direct-coordinate} therefore yields 
\[
(n+2)\overline{B}_{p}(R)=Q_{p}\left(3R+\frac{1}{2}\Ric(R)\owedge g\right)-\frac{1}{2}\bigl(n\Ric(R)-\Scal(R)g\bigr)^{[p]}.
\]
Equation \eqref{eq:kulkarni-q} and $g^{[p]}=p\Id$ reduce this to
\eqref{eq:average}. 
\end{proof}

\section{The curvature estimate and the homology sphere theorem}

\label{sec:proof}
\begin{proof}[Proof of Theorem~\ref{thm:algebraic}]
First, $\Ric_{p}(E)\geq0$ implies $\Ric_{j}(E)\geq0$ for $p\leq j\leq n-1$.
Indeed, for an orthonormal frame $v,e_{1},\ldots,e_{j}$, summing
over its $p$-element subsets gives 
\[
0\leq\sum_{\substack{A\subset\{1,\ldots,j\}\\
|A|=p
}
}\sum_{i\in A}E(v,e_{i},v,e_{i})=\binom{j-1}{p-1}\sum^{j}_{i=1}E(v,e_{i},v,e_{i}).
\]
Division by the positive binomial coefficient proves the claim, and
the same argument preserves strict positivity. Since $p\leq n/2$,
it applies in particular to $j=n-p$.

On every $(p+1)$-plane $H$, the hypothesis yields $\Ric(E|_{H})\geq0$.
Hodge duality within $H$ identifies this operator with $Q_{p}(E|_{H})$,
by \eqref{eq:q-ricci} and \eqref{eq:q-hodge}. Extending by zero
and summing therefore gives 
\begin{equation}
B^{(p)}_{\mathcal{B}}(E)\geq0 \text{ and } \overline{B}_{p}(E)\geq0.\label{eq:positive-average}
\end{equation}
The subset calculation gives the same conclusion in degree $n-p$.
The two averaged identities are consequently 
\begin{equation}\label{eq:first-positive}
    3Q_{p}(E)-p\Ric(E)^{[p]}+p\Scal(E)\Id  =(n+2)\overline{B}_{p}(E)\geq0,
\end{equation}
and 
\begin{equation}\label{eq:second-positive}
3Q_{p}(E)+(n-p)\Ric(E)^{[p]}  =(n+2)*^{-1}\overline{B}_{n-p}(E)*\geq0.
\end{equation}
The cancellation identity \eqref{eq:positive-representation} thus
has a nonnegative right-hand side. Dividing by $3n$ proves \eqref{eq:algebraic}.

If $\Ric_{p}(E)>0$, each $Q_{p}(E|_{H})$ is positive definite on
$\Lambda^{p}H^{*}$. For any frame and any nonzero $p$-form $\omega$,
choose a nonzero coefficient of $\omega$ in the basis $e^{A}$ and
a coordinate $(p+1)$-plane containing those directions. The restriction
of $\omega$ to this plane is nonzero, so its contribution to $\langle B^{(p)}_{\mathcal{B}}(E)\omega,\omega\rangle$
is positive. All other contributions are nonnegative. Thus $B^{(p)}_{\mathcal{B}}(E)$
and its average are positive definite, making the right-hand side
of \eqref{eq:positive-representation} positive definite. This proves
the assertion when the pinching is strict at a point.

We now prove \eqref{eq:main-estimate} in all degrees. Set 
\[
E=R-\frac{1}{2}\Ric_{2}(R)I.
\]
Then $\Ric_{2}(E)=0$, and the subset calculation shows that $\Ric_{p}(E)\geq0$
for $2\leq p\leq n-1$. The restricted Ricci operators therefore imply
that $B^{(p)}_{\mathcal{B}}(E)\geq0$ for $2\leq p\leq n-1$. In degree
1, summing the operators of the coordinate two-planes gives 
\[
B^{(1)}_{\mathcal{B}}(E)=\Ric(E)_{\mathrm{d}}\geq0,
\]
since $\Ric_{n-1}(E)\geq0$. Hence $\overline{B}_{p}(E)\geq0$ for
every $1\leq p\leq n-1$. The identity in \eqref{eq:positive-representation}
is valid in every such degree, so it yields 
\[
Q_{p}(E)\geq-\frac{p(n-p)}{3n}\Scal(E)\Id.
\]
Using $\Scal(E)=\Scal(R)-\frac{n(n-1)}{2}\Ric_{2}(R)$ and $Q_{p}(I)=p(n-p)\Id$,
we obtain 
\begin{align*}
Q_{p}(R) & =Q_{p}(E)+\frac{p(n-p)}{2}\Ric_{2}(R)\Id\\
 & \geq\frac{p(n-p)}{3n}\left(\frac{3n}{2}\Ric_{2}(R)-\Scal(E)\right)\Id\\
 & =\frac{p(n-p)}{3n}\left(\frac{n(n+2)}{2}\Ric_{2}(R)-\Scal(R)\right)\Id.
\end{align*}
This proves \eqref{eq:main-estimate} for every degree $p=1,...,n-1$. 
\end{proof}

\begin{proof}[Proof of Theorem~\ref{thm:intermediate}]
At each point, set 
\[
E=R-\frac{\Scal(R)}{n(n+2)}I.
\]
The hypothesis and the shift identity for $\Ric_{p}$ imply that 
\[
\Ric_{p}(E)=\Ric_{p}(R)-\frac{p}{n(n+2)}\Scal(R)\geq0.
\]
The subset calculation then yields $\Ric_{j}(E)\geq0$ for $p\leq j\leq n-1$.
Moreover, 
\[
\Scal(E)=\Scal(R)-\frac{n-1}{n+2}\Scal(R)=\frac{3}{n+2}\Scal(R).
\]
For $p\leq j\leq\lfloor\frac{n}{2}\rfloor$, Theorem~\ref{thm:algebraic}
and \eqref{eq:q-constant} therefore show that 
\begin{align*}
Q_{j}(R) & =Q_{j}(E)+\frac{j(n-j)}{n(n+2)}\Scal(R)\Id\\
 & \geq-\frac{j(n-j)}{3n}\Scal(E)\Id+\frac{j(n-j)}{n(n+2)}\Scal(R)\Id\\
 & =0.
\end{align*}
Hodge duality in each tangent space gives the same conclusion in the
complementary degrees, so $Q_{j}(R)\geq0$ for $p\leq j\leq n-p$.
The Bochner identity \eqref{eq:bochner} makes every harmonic form
in this range parallel.

At a point where \eqref{eq:intermediate-pinching} is strict, we have
$\Ric_{p}(E)>0$. The subset calculation preserves positivity, so
the strict part of Theorem~\ref{thm:algebraic}, followed by the
same scalar cancellation and Hodge duality, makes $Q_{j}(R)$ positive
definite there for $p\leq j\leq n-p$. A parallel harmonic form must
therefore vanish there, and hence on all of $M$. Hodge theory proves
the asserted vanishing. Only Hodge duality in each tangent space was
used, so orientability is not required. 
\end{proof}

\begin{proof}[Proof of Corollary~\ref{cor:second-ricci-strict}]
The hypothesis and \eqref{eq:main-estimate} give $Q_{j}(R)\geq0$
in every intermediate degree, with positive definiteness at a point
where the pinching inequality is strict. The Bochner identity \eqref{eq:bochner}
first makes every harmonic form parallel, and positivity at that point
forces it to vanish. Hodge theory gives the asserted vanishing, without
an orientability assumption. 
\end{proof}

\begin{proof}[Proof of Theorem~\ref{thm:main}]
Since $\frac{1}{2}\Ric_{2}(R)\geq\min\sec(R)$, the sectional pinching
hypothesis implies the condition of Corollary~\ref{cor:second-ricci-strict}.
Thus $H^{p}(M;\R)=0$ for every degree $=1,...,n-1$, whether or not
$M$ is orientable.

A closed manifold has the homotopy type of a finite CW complex. Extension
of coefficients therefore shows that, for $1\leq p\leq n-1$, 
\[
H^{p}(M;\Q)\otimes_{\Q}\R\cong H^{p}(M;\R)=0.
\]
Consequently $H^{p}(M;\Q)=0$, and the universal coefficient theorem
over $\Q$ gives $H_{p}(M;\Q)=0$ in the same degrees. If $M$ is
connected and orientable, its remaining rational homology groups are
$H_{0}(M;\Q)\cong H_{n}(M;\Q)\cong\Q$, with all groups above degree
$n$ equal to zero. Thus $M$ is a rational homology sphere. When
$n$ is odd, the Euler characteristic argument in Section~\ref{sec:prelim}
gives orientability, so the same conclusion holds. 
\end{proof}

\begin{remark}[Sharpness]\label{rem:sharpness} In every even dimension
$n=2m\geq4$, the Fubini-Study metric on $\CP^{m}$ satisfies $\min\sec=1$,
$\Scal=n(n+2)$, and $b_{2}=1$. Thus the strict inequality in Theorem~\ref{thm:main}
cannot be replaced by a weak inequality, and the coefficient $\frac{1}{n(n+2)}$
cannot be decreased in a rational homology sphere theorem of this
form in any even dimension. The algebraic estimate \eqref{eq:algebraic}
is also attained for each even $p$, $2\leq p\leq n-2$: take $E=R_{\mathrm{FS}}-I$
and let $\omega$ be a nonzero power of the parallel Kähler form.
Since $Q_{p}(R_{\mathrm{FS}})\omega=0$, we have 
\[
\Scal(E)=3n
\]
and 
\[
Q_{p}(E)\omega=-p(n-p)\omega=-\frac{p(n-p)}{3n}\Scal(E)\omega.
\]
These examples do not assert sharpness for odd degrees or odd ambient
dimensions. \end{remark}

\begin{remark} The strict inequality \eqref{eq:yau} implies $\sec>0$,
since $\Scal\geq n(n-1)\min\sec$. In even dimensions, cohomology
vanishing in degrees $1,...,n-1$ does not determine top-degree cohomology:
a round $\mathbb{RP}^{2m}$ satisfies \eqref{eq:yau} but has zero
top rational homology. Round odd-dimensional real projective spaces
satisfy the same curvature condition and have integral torsion, so
an integral homology sphere conclusion would require additional hypotheses.
In dimensions 2 and 3, the real cohomology conclusion follows directly
from positive Ricci curvature and Hodge duality in each tangent space.
\end{remark}

\begin{remark}[Comparison with spectral conditions]\label{rem:spectral-comparison}
The following algebraic examples in dimension 7 distinguish our hypotheses
from the homology sphere criteria in \cite{PW,NPW,LiCOSK,LPCOSK}.
For the Fubini-Study tensor with $\min\sec=1$ on an even-dimensional
subspace $H\subset\R^{7}$ equipped with an orthogonal complex structure,
the tensor $R_{\mathrm{FS}}-I|_{H}$ has nonnegative sectional curvature
and Ricci tensor $3g|_{H}$. We extend this tensor by zero to $\R^{7}$.
The spectra used below follow from the Fubini-Study curvature formula;
see also \cite[Appendix~A.5]{LPCOSK}.

Take $\dim H=6$ and set 
\[
R=I+\frac{13}{12}(R_{\mathrm{FS}}-I|_{H}).
\]
Then $\min\sec(R)=1$ and $\Scal(R)=\frac{123}{2}<63$, so \eqref{eq:yau}
holds. The curvature operator has eigenvalue $-\frac{1}{12}$ with
multiplicity six and hence is not $4$-nonnegative. On the 8-dimensional
space of trace-free Hermitian symmetric two-tensors supported on $H$,
the curvature operator of the second kind has eigenvalue $-\frac{9}{4}$.
The modification in \cite[formula~(1.5)]{LPCOSK} changes this eigenvalue
to 
\[
-\frac{9}{4}+\frac{4}{11}\left(\frac{37}{4}-\frac{123}{14}\right)+\frac{123}{63}=-\frac{17}{132}.
\]
Thus neither operator of the second kind is $\frac{9}{2}$-nonnegative.

The example also fails the enlarged homology criterion in \cite[Theorem~1.6]{LiCOSK}:
in dimension 7 its limiting lower bound for the average of the lowest
$\frac{9}{2}$ eigenvalues is $-\frac{\Scal}{189}$, which exceeds
$-\frac{9}{4}$.

Conversely, let $L=\operatorname{diag}(1,1,1,1,1,1,5)$ and $R=L\owedge g$.
The curvature operator is positive, and the quadratic form of the
curvature operator of the second kind on a trace-free symmetric two-tensor
$h$ is $2\tr(Lh^{2})>0$ for $h\ne0$. The corresponding quadratic
form for the modified operator in \cite{LPCOSK} is 
\[
\frac{42}{11}\tr(Lh^{2})+\frac{4}{3}|h|^{2}>0.
\]
Nevertheless, $\Scal(R)=132$ and $\min\sec(R)=2<\frac{\Scal(R)}{7\cdot9}=\frac{44}{21}$,
so \eqref{eq:yau} fails. Moreover, $\Ric_{p}(R)=2p$ for $1\leq p\leq3$,
so \eqref{eq:intermediate-pinching} fails as well. The first example
satisfies all these lower intermediate Ricci bounds. These examples
show that the various curvature conditions do not imply each other,
without asserting it in every dimension or for every fixed parameter
in Li's family. \end{remark}

\section{Weak pinching and holonomy}

\label{sec:odd}

We combine the fact that harmonic forms are parallel with the holonomy
classification to prove Theorems~\ref{thm:intermediate-rigidity}
and \ref{thm:odd-weak}. We first establish the curvature consequences
needed to exclude products and to handle the Kähler and quaternionic
Kähler cases.

\begin{lemma}\label{lem:intermediate-ricci-consequences} Let $R$
be an algebraic curvature tensor on an $n$-dimensional Euclidean
space, $n\geq4$, and let $1\leq p\leq n-2$. If 
\[
\Ric_{p}(R)\geq\frac{p\Scal(R)}{n(n+2)},
\]
then 
\[
\Ric(R)\geq\frac{n-1}{n(n+2)}\Scal(R)g
\]
and 
\[
\Scal(R)\geq0.
\]
Moreover, $\Scal(R)=0$ if and only if $R=0$. If a closed connected
nonflat Riemannian manifold satisfies \eqref{eq:intermediate-pinching}
for $1\leq p\leq\lfloor\frac{n}{2}\rfloor$, its universal Riemannian
cover has irreducible holonomy. \end{lemma}
\begin{proof}
Fix a unit vector $v$ and an orthonormal basis $e_{1},\ldots,e_{n-1}$
of $v^{\perp}$. Summing over its $p$-element subsets gives 
\begin{align*}
\binom{n-2}{p-1}\Ric(R)(v,v) & =\sum_{\substack{A\subset\{1,\ldots,n-1\}\\
|A|=p
}
}\sum_{i\in A}R(v,e_{i},v,e_{i}) \geq\binom{n-1}{p}\frac{p\Scal(R)}{n(n+2)}.
\end{align*}
Since $\frac{\binom{n-1}{p}}{\binom{n-2}{p-1}}=\frac{n-1}{p}$, we
obtain 
\[
\Ric(R)(v,v)\geq\frac{n-1}{n(n+2)}\Scal(R).
\]
Taking the trace yields 
\[
\Scal(R)\geq\frac{n-1}{n+2}\Scal(R),
\]
so $\Scal(R)\geq0$.

If $\Scal(R)=0$, the Ricci tensor is nonnegative and has zero trace,
and therefore vanishes. Fix a unit vector $v$ and an orthonormal
basis of $v^{\perp}$. Every $p$-term sum of the sectional curvatures
through $v$ in this basis is nonnegative, and the sum of all such
sums is 
\[
\binom{n-2}{p-1}\Ric(R)(v,v)=0.
\]
Hence each $p$-term sum is zero. Comparing two subsets that differ
in one index, which is possible since $p\leq n-2$, shows that these
sectional curvatures are all equal. Their total is zero, so each vanishes.
Since the frame was arbitrary, every sectional curvature is zero and
$R=0$. The converse is immediate.

For the last assertion, nonflatness now gives a point of positive
scalar curvature. The universal Riemannian cover is complete. If it
splits as a nontrivial Riemannian product, one factor would have dimension
at most $\frac{n}{2}$, and its orthogonal complement would have dimension
at least $p$. Choose $v$ tangent to this factor and orthonormal
$e_{1},\ldots,e_{p}$ tangent to the complement. All $p$ mixed sectional
curvatures vanish, contradicting \eqref{eq:intermediate-pinching}
at a point of positive scalar curvature. The de Rham decomposition
theorem therefore makes the holonomy representation irreducible. 
\end{proof}

\begin{lemma}\label{lem:kahler-endpoint} Let $R$ be an algebraic
Kähler curvature tensor on a Hermitian space $V$ of real dimension
$n\geq4$. Suppose that 
\[
\Ric_{n-2}(R)\geq\frac{n-2}{n(n+2)}\Scal(R).
\]
It follows that $R$ is the curvature tensor of a complex space form
with nonnegative holomorphic sectional curvature $\frac{4\Scal(R)}{n(n+2)}$.
\end{lemma}
\begin{proof}
Write $J$ for the orthogonal complex structure. We use the Kähler
symmetries 
\[
R(JX,JY,Z,W)=R(X,Y,Z,W)=R(X,Y,JZ,JW).
\]
Put $k=\frac{\Scal(R)}{n(n+2)}$. Lemma~\ref{lem:intermediate-ricci-consequences}
shows that $k\geq0$. For a unit vector $X$, let $H(X)=R(X,JX,X,JX)$.
Equation \eqref{eq:sphere-quartic} and the Kähler symmetries yield
\begin{equation}
\int_{S(V)}H(X)\,d\sigma(X)=\frac{4\Scal(R)}{n(n+2)}=4k.\label{eq:holomorphic-average}
\end{equation}
Indeed we have the three standard contractions
\begin{align*}
\sum_{a,b}R(e_{a},Je_{a},e_{b},Je_{b}) & =\sum_{a,b}R(e_{a},e_{b},e_{a},e_{b})+\sum_{a,b}R(e_{a},Je_{b},e_{a},Je_{b}) =2\Scal(R),
\end{align*}
\[
\sum_{a,b}R(e_{a},Je_{b},e_{a},Je_{b})=\Scal(R),
\]
and 
\[
\sum_{a,b}R(e_{a},Je_{b},e_{b},Je_{a})=\sum_{a,b}R(e_{a},Je_{b},e_{a},Je_{b})=\Scal(R).
\]
Together they give the usual holomorphic sectional curvature average;
see also \cite{Berger1966,Magnusson}.

Complete $X,JX$ to an orthonormal basis $X,JX,e_{3},\ldots,e_{n}$.
The pinching assumption shows that 
\[
\Ric(R)(X,X)-H(X)=\sum^{n}_{a=3}R(X,e_{a},X,e_{a})\geq(n-2)k.
\]
Equations \eqref{eq:sphere-quadratic} and \eqref{eq:holomorphic-average}
give 
\[
\int_{S(V)}\bigl(\Ric(R)(X,X)-H(X)-(n-2)k\bigr)\,d\sigma(X)=0.
\]
The integrand is continuous and nonnegative, so 
\begin{equation}
\Ric(R)(X,X)-H(X)=(n-2)k\label{eq:orthogonal-ricci-equality}
\end{equation}
for every unit $X$.

For any unit $Y\perp X$, apply the hypothesis to an orthonormal basis
of $\{X,Y\}^{\perp}$. It gives 
\[
\Ric(R)(X,X)-R(X,Y,X,Y)\geq(n-2)k.
\]
Together with \eqref{eq:orthogonal-ricci-equality}, this shows that
\[
R(X,Y,X,Y)\leq H(X).
\]
Thus the holomorphic plane maximizes sectional curvature among planes
containing $X$. Fix a unit vector $Y\perp\operatorname{span}\{X,JX\}$
and apply the same observation with first vector $JX$. Differentiating
at the maximizing direction $X$ gives 
\begin{align*}
0 & =\left.\frac{d}{dt}\right|_{t=0}R(JX,\cos t\,X+\sin t\,Y,JX,\cos t\,X+\sin t\,Y)\\
 & =2R(JX,Y,JX,X).
\end{align*}
The Kähler symmetries then give 
\begin{align*}
\left.\frac{d}{dt}\right|_{t=0}H(\cos t\,X+\sin t\,Y) & =2R(Y,JX,X,JX)+2R(X,JY,X,JX)\\
 & =4R(X,JX,Y,JX)\\
 & =4R(JX,Y,JX,X) \\
 & =0.
\end{align*}
The derivative in the direction $JX$ also vanishes, since $H(\cos t\,X+\sin t\,JX)=H(X)$.
Thus every tangential derivative of $H$ on $S(V)$ is zero. Its average
in \eqref{eq:holomorphic-average} gives $H=4k$. A Kähler curvature
tensor is determined by its holomorphic sectional curvatures, so polarization
identifies $R$ with the curvature tensor of a complex space form
with holomorphic sectional curvature $4k$. 
\end{proof}

\begin{remark} The hypothesis of Lemma~\ref{lem:kahler-endpoint}
follows from the analogous bound on $\Ric_{p}(R)$ for any $1\leq p\leq n-2$,
by summing over subsets. The index $n-2$ cannot be replaced by $n-1$:
every positive Kähler-Einstein metric satisfies 
\[
\Ric=\frac{\Scal}{n}g>\frac{n-1}{n(n+2)}\Scal\,g,
\]
including products of positive Kähler-Einstein manifolds with the
same Einstein constant. \end{remark}

\begin{lemma}\label{lem:quaternionic-bound} On a quaternionic Kähler
manifold of real dimension $n=4m\geq8$, 
\begin{equation}
\min\sec\leq\frac{\Scal}{n(n+8)}.\label{eq:quaternionic-bound}
\end{equation}
If, for some $1\leq p\leq n-4$, the metric satisfies \eqref{eq:intermediate-pinching},
then it is flat. \end{lemma}
\begin{proof}
Choose orthogonal complex structures $J_{1},J_{2},J_{3}$ spanning
the quaternionic structure at a point, with $J_{1}J_{2}=J_{3}=-J_{2}J_{1}$.
Since $n\geq8$, there are unit vectors $X,Y$ such that 
\[
Y\perp\operatorname{span}\{X,J_{1}X,J_{2}X,J_{3}X\}.
\]
The four vectors $Y,J_{1}Y,J_{2}Y,J_{3}Y$ are orthonormal and perpendicular
to $X$. The quaternionic curvature identity \cite[Theorems~1.2 and~1.3(2)]{KLZ}
gives 
\begin{align*}
4\min\sec  \leq R(X,Y,X,Y)+\sum^{3}_{a=1}R(X,J_{a}Y,X,J_{a}Y) =\frac{4\Scal}{n(n+8)}.
\end{align*}
This proves \eqref{eq:quaternionic-bound}. To check the normalization,
the Einstein identity in \cite[Theorem~1.2]{KLZ}, written in real
dimension $n=4m$, is 
\[
\Ric=4(m+2)\delta g=(n+8)\delta g.
\]
Taking the trace gives $\Scal=n(n+8)\delta$, so the constant $\delta$
is $\frac{\Scal}{n(n+8)}$ in our notation.

Now assume \eqref{eq:intermediate-pinching} for $1\leq p\leq n-4$.
The other quaternionic curvature identity \cite[Theorem~1.3(1)]{KLZ}
gives 
\[
\sum^{3}_{a=1}R(X,J_{a}X,X,J_{a}X)=\frac{12\Scal}{n(n+8)}.
\]
Choose an orthonormal basis $e_{1},\ldots,e_{n-4}$ of $\operatorname{span}\{X,J_{1}X,J_{2}X,J_{3}X\}^{\perp}$.
Subtracting the displayed identity from $\Ric(X,X)=\frac{\Scal}{n}$
yields 
\[
\sum^{n-4}_{i=1}R(X,e_{i},X,e_{i})=\frac{n-4}{n(n+8)}\Scal.
\]
Since $p\leq n-4$, summing the assumed $p$-term inequalities over
subsets bounds the left-hand side below by $\frac{n-4}{n(n+2)}\Scal$.
Consequently, 
\[
\frac{\Scal}{n(n+8)}\geq\frac{\Scal}{n(n+2)}.
\]
Lemma~\ref{lem:intermediate-ricci-consequences} gives $\Scal\geq0$,
so this inequality forces $\Scal=0$. The same lemma then gives flatness. 
\end{proof}

\begin{proof}[Proof of Theorem~\ref{thm:intermediate-rigidity}]
By Lemma~\ref{lem:intermediate-ricci-consequences}, the universal
Riemannian cover has irreducible holonomy, and the Ricci tensor is
positive at every nonflat point. In particular, the metric is not
Ricci flat.

Suppose first that the metric is not locally symmetric. The Berger-Simons
classification \cite{Berger,Simons} leaves the following possibilities
for restricted holonomy: 
\begin{align*}
 & \mathrm{SO}(n),\\
 & \mathrm{U}(m),\text{ for }n=2m,\\
 & \mathrm{Sp}(m)\mathrm{Sp}(1),\text{ for }n=4m\geq8.
\end{align*}
The groups $\mathrm{SU}(m)$, $\mathrm{Sp}(m)$, $\mathrm{G}_{2}$,
and $\mathrm{Spin}(7)$ force Ricci flatness and are therefore excluded;
see \cite[10.64]{Besse}. In the quaternionic dimensions, $p\leq\frac{n}{2}\leq n-4$,
so Lemma~\ref{lem:quaternionic-bound} excludes $\mathrm{Sp}(m)\mathrm{Sp}(1)$
as well.

If the restricted holonomy is $\mathrm{U}(m)$, the simply connected
universal cover carries a parallel orthogonal complex structure. Since
$p\leq n-2$, summing over subsets gives the hypothesis of Lemma~\ref{lem:kahler-endpoint}.
Its holomorphic sectional curvature is constant at each point, with
value $\frac{4\Scal}{n(n+2)}$. In particular, $\Ric=\frac{\Scal}{n}g$.
The contracted Bianchi identity gives 
\[
\frac{1}{n}\,d\Scal=\operatorname{div}\Ric=\frac{1}{2}\,d\Scal.
\]
Since $n\geq4$, scalar curvature is constant. It is positive because
the metric is nonflat. The metric is therefore locally a scaled Fubini-Study
metric, contrary to the assumption that it is not locally symmetric.
Thus the remaining nonsymmetric restricted holonomy is $\mathrm{SO}(n)$.

A parallel form is fixed by restricted holonomy. For $1\leq j\leq n-1$,
take any orthonormal $j$-frame and rotate through $\pi$ in the plane
spanned by its first vector and a unit vector perpendicular to the
frame. This rotation lies in $\mathrm{SO}(n)$ and changes the sign
of the form's evaluation on the frame. Invariance therefore forces
that evaluation to vanish. Since the frame was arbitrary, the form
is zero. Theorem~\ref{thm:intermediate} makes every harmonic form
of degree $p\leq j\leq n-p$ parallel. Thus $b_{j}(M)=0$ in this
range.

If some $b_{j}(M)\ne0$ in this range, the preceding argument forces
local symmetry. The complete simply connected universal cover is then
an irreducible symmetric space. Its Ricci tensor is positive at a
point and is parallel, so it is positive everywhere. Thus the symmetric
space is of compact type and is compact; see \cite[10.80]{Besse}.
This proves the final assertion. 
\end{proof}

\begin{proof}[Proof of Theorem~\ref{thm:odd-weak}]
Lemma~\ref{lem:intermediate-ricci-consequences} with $p=2$ shows
that the Ricci curvature is nonnegative and positive at a nonflat
point. The Bochner formula makes every harmonic one-form parallel,
and positivity at that point forces it to vanish. This argument also
applies on the orientation cover. Hodge duality there and injectivity
of pullback give $b_{1}(M)=b_{n-1}(M)=0$. If all the Betti numbers
in degrees $2,\ldots,n-2$ vanish, we obtain the first alternative.
Otherwise, Theorem~\ref{thm:intermediate-rigidity} with $p=2$ shows
that the universal Riemannian cover $N$ is a compact irreducible
symmetric space. Its nonflat symmetric metric has positive constant
scalar curvature.

The rank of $N$ must be one. Otherwise, choose a unit vector $X$
on a restricted-root hyperplane in the tangent space of a maximal
flat. Such a vector exists because the flat has dimension at least
two. Choose a unit vector $e_{1}$ tangent to the flat and perpendicular
to $X$, and a unit vector $e_{2}$ in the tangent root space corresponding
to that hyperplane. The root-space decomposition makes $e_{2}$ perpendicular
to the flat. Both $e_{1}$ and $e_{2}$ commute with $X$ in the symmetric-space
Lie algebra, the latter because the corresponding root vanishes on
$X$; see \cite[Section~4.1]{DGM}. Commuting tangent vectors span
flat planes, so this orthonormal triple satisfies 
\[
\Ric_{2}(R)\leq\sec(X,e_{1})+\sec(X,e_{2})=0.
\]
This contradicts \eqref{eq:second-ricci-pinching}, since $\Scal>0$.

The simply connected compact rank-one symmetric spaces are round spheres,
complex and quaternionic projective spaces, and the Cayley plane;
see \cite[10.80-10.81]{Besse}. Quaternionic projective spaces of
real dimension at least eight are excluded by Lemma~\ref{lem:quaternionic-bound}
with $p=2$; the quaternionic projective line is a round four-sphere.
For the Cayley plane normalized by $\min\sec=1$ and any unit vector
$X$, the Jacobi operator on $X^{\perp}$ has eigenvalues $1$ and
$4$ with multiplicities eight and seven, respectively; see \cite[Section~3]{Nikolayevsky}.
Thus $\Scal=16(8+4\cdot7)=576$. Choosing orthonormal $e_{1},e_{2}$
in the eigenvalue-one space gives 
\[
\Ric_{2}(R)\leq\sec(X,e_{1})+\sec(X,e_{2})=2<\frac{2\cdot576}{16\cdot18}=4,
\]
so \eqref{eq:second-ricci-pinching} fails.

The remaining models are round spheres and scaled Fubini-Study complex
projective spaces. A compact spherical cover is finite, and pullback
in real cohomology is injective. In case of the sphere it therefore forces
all intermediate Betti numbers of $M$ to vanish, contrary to our
assumption. Thus $N$ is a scaled Fubini-Study complex projective
space, giving the complex-projective quotient in the second alternative.
If $n$ is odd, this alternative is impossible. The first alternative
and the observation in Section~\ref{sec:prelim} then give the rational
homology sphere conclusion. 
\end{proof}

\begin{remark} Flat tori show that the nonflat assumption is necessary,
and Fubini-Study complex projective spaces account for the exception
in even dimensions. The quaternionic exclusion is algebraic; it does
not require nonnegativity of the curvature operator or an Einstein
rigidity theorem. In the closed setting and under the stronger sectional
assumption $\sec\geq\frac{\Scal}{n(n+2)}$, the Einstein property
of quaternionic Kähler metrics in dimension at least eight and the
weak pinching theorem of Xu and Gu \cite[Theorem~1.1(ii)]{XG} offer
an alternative exclusion of the nonsymmetric $\mathrm{Sp}(m)\mathrm{Sp}(1)$
case. \end{remark}

\begin{remark}[An improvement in dimension five] \label{rem:odd-low-dimensional}
A closed connected 5-manifold with a nonflat metric satisfying 
\[
\sec\geq\frac{\Scal}{38}
\]
at every point is a rational homology sphere. This coefficient is
smaller than Yau's coefficient $\frac{1}{35}$ and follows from Ge's
two-form estimate \cite[Theorem~1.3]{Ge}. Indeed, tracing the sectional
bound gives $\Scal\geq0$, and $\Scal=0$ forces every sectional curvature
to vanish. Thus 
\[
\Ric\geq\frac{2}{19}\Scal\,g\geq0,
\]
with positive definiteness at a nonflat point. The Bochner formula
eliminates harmonic 1-forms on $M$ and on its orientation cover.
On that cover, Ge's theorem gives $b_{2}=0$, and Hodge duality gives
$b_{3}=b_{4}=0$. Injectivity of pullback and the Euler characteristic
argument in Section~\ref{sec:prelim} prove the claim.

This improvement and the analogy with quarter pinching in \cite{BS, PT, NW} suggest that
$\frac{1}{n(n+2)}$ can be decreased in every odd dimension $n\geq5$.
We leave this question open. \end{remark}

\section{Eigenvalue estimates and spectral pinching}

\label{sec:spectral-pinching}

The estimate in Theorem~\ref{thm:algebraic} gives lower bounds for
the Hodge Laplacian and a vanishing criterion involving a scalar Schrödinger
operator. By the definition in Section~\ref{sec:prelim}, $\Ric_{2}(R)$
is the minimum of $\sec(v,e_{1})+\sec(v,e_{2})$ over orthonormal
triples. Put 
\[
D_{2}=\frac{n(n+2)}{2}\Ric_{2}(R)-\Scal.
\]
The function $D_{2}$ is continuous. The same notation applies to
an algebraic curvature tensor on a Euclidean space.

Equation \eqref{eq:main-estimate} gives, in every intermediate degree,
\[
Q_{j}(R)\geq\frac{j(n-j)}{3n}D_{2}\Id.
\]

\begin{theorem}\label{thm:hodge-eigenvalue} Let $(M^{n},g)$ be
a closed connected Riemannian manifold, $n\geq4$, with $D_{2}>0$
everywhere. For $1\leq j\leq n-1$, let $\lambda_{1,j}$ denote the
least positive eigenvalue of the Hodge Laplacian on $j$-forms. Then
\begin{equation}
\lambda_{1,j}\geq\frac{\min\bigl\{ j(n-j+1),(j+1)(n-j)\bigr\}}{3n}\min_{M}D_{2}.\label{eq:hodge-eigenvalue}
\end{equation}
There are no nonzero harmonic forms in these degrees. The estimate
is attained on round spheres. \end{theorem}
\begin{proof}
Theorem~\ref{thm:algebraic} and the Bochner identity \eqref{eq:bochner}
eliminate harmonic $j$-forms. We give the usual argument for the
eigenvalue estimate; see \cite{GM}.

For a $j$-form $\omega$, the orthogonal projections of $\nabla\omega$
onto its alternation and contraction parts have squared norms $\frac{1}{j+1}|d\omega|^{2}$
and $\frac{1}{n-j+1}|d^{*}\omega|^{2}$, respectively. Consequently,
\[
|\nabla\omega|^{2}\geq\frac{1}{j+1}|d\omega|^{2}+\frac{1}{n-j+1}|d^{*}\omega|^{2}.
\]
Combining this inequality with the integrated Weitzenböck formula
gives 
\begin{align*}
\frac{j}{j+1}\|d\omega\|^{2}_{L^{2}}+\frac{n-j}{n-j+1}\|d^{*}\omega\|^{2}_{L^{2}}  \geq\int_{M}\langle Q_{j}(R)\omega,\omega\rangle\,dV
  \geq\frac{j(n-j)}{3n}(\min_{M}D_{2})\|\omega\|^{2}_{L^{2}}.
\end{align*}
If $\omega$ is an eigenvector with eigenvalue $\lambda>0$, then
$\|d\omega\|^{2}_{L^{2}}+\|d^{*}\omega\|^{2}_{L^{2}}=\lambda\|\omega\|^{2}_{L^{2}}$.
Thus 
\[
\max\left\{ \frac{j}{j+1},\frac{n-j}{n-j+1}\right\} \lambda\geq\frac{j(n-j)}{3n}\min_{M}D_{2},
\]
which is \eqref{eq:hodge-eigenvalue}.

On the unit round sphere, $\Ric_{2}(R)=2$ and $\Scal=n(n-1)$, so
$D_{2}=3n$. Constant forms on $\R^{n+1}$ produce eigenvectors attaining
both values in the minimum. If $\alpha$ is a nonzero constant $r$-form,
$1\leq r\leq n$, let $\phi$ be its tangential restriction to $S^{n}$
and let $\psi$ be the tangential restriction of $\iota_{x}\alpha$,
where $x$ is the position vector. The Gauss formula gives 
\[
\nabla_{X}\phi=-X^{\flat}\wedge\psi
\]
and 
\[
\nabla_{X}\psi=\iota_{X}\phi,
\]
where $X^{\flat}=g(X,\cdot)$. Consequently, 
\begin{align*}
d\phi & =0, \\
 d^{*}\phi & =(n-r+1)\psi,\\
d\psi & =r\phi, \\
 d^{*}\psi & =0.
\end{align*}
Taking $r=j$ gives a nonzero exact $j$-eigenform with eigenvalue
$j(n-j+1)$. Taking $r=j+1$ gives a nonzero coexact $j$-eigenform
with eigenvalue $(j+1)(n-j)$. This proves sharpness. Scaling gives
equality on every round sphere. 
\end{proof}

The next result allows $D_{2}$ to have either sign. Write $\Delta=d^{*}d$
for the nonnegative Laplacian on functions. For a continuous real-valued
potential $V$, define 
\[
\lambda_{0}(\Delta+V)=\inf_{\substack{u\in W^{1,2}(M;\R)\\
u\ne0
}
}\frac{{\displaystyle \int_{M}\bigl(|du|^{2}+Vu^{2}\bigr)\,dV}}{{\displaystyle \int_{M}u^{2}\,dV}}.
\]
Here $W^{1,2}(M;\R)$ is the Sobolev space of real-valued functions
with square-integrable weak differential.

\begin{theorem}\label{thm:spectral-pinching} Let $(M^{n},g)$ be
a closed connected Riemannian manifold, $n\geq4$, and set 
\[
C_{n}=\begin{cases}
{\displaystyle \frac{n^{2}}{12(n+2)},} & n\text{ even},\\[6pt]
{\displaystyle \frac{(n-1)(n+1)^{2}}{12n(n+3)},} & n\text{ odd}.
\end{cases}
\]
If 
\begin{equation}
\lambda_{0}(\Delta+C_{n}D_{2})>0,\label{eq:spectral-pinching}
\end{equation}
then $H^{j}(M;\R)=0$ for $1\leq j\leq n-1$. If $n$ is odd or $M$
is orientable, it is a rational homology sphere. \end{theorem}
\begin{proof}
We first verify that 
\[
C_{n}=\max_{1\leq j\leq n-1}\frac{j(n-j)}{3n}\frac{\max\{j,n-j\}}{\max\{j,n-j\}+1}.
\]
By symmetry, it suffices to consider $1\leq j\leq\lfloor\frac{n}{2}\rfloor$.
The expression then equals $\frac{j(n-j)^{2}}{3n(n-j+1)}$ and increases
with $j$. Indeed, for $j<\lfloor\frac{n}{2}\rfloor$, put $q=n-j\geq j+2$.
Then 
\begin{align*}
\frac{(j+1)(q-1)^{2}}{q}-\frac{jq^{2}}{q+1}  =\frac{(q-j-1)q^{2}-(j+1)(q-1)}{q(q+1)}
  \geq\frac{2q-1}{q(q+1)}>0.
\end{align*}
Thus the maximum occurs at the middle degree or degrees, giving the
stated values of $C_{n}$.

Fix $1\leq j\leq n-1$ and put $q=\max\{j,n-j\}$. Since 
\[
0<\frac{j(n-j)q}{3n(q+1)}\leq C_{n},
\]
the quadratic form of $\Delta+\frac{j(n-j)q}{3n(q+1)}D_{2}$ is a
convex combination of those of $\Delta$ and $\Delta+C_{n}D_{2}$,
with positive weight on the latter. The first form is nonnegative,
and the second is positive by \eqref{eq:spectral-pinching}. Hence
\[
\int_{M}\left(|du|^{2}+\frac{j(n-j)q}{3n(q+1)}D_{2}u^{2}\right)\,dV>0
\]
for every nonzero \(u\in W^{1,2}(M;\R) \). 

Suppose that $\omega$ is a nonzero harmonic $j$-form. Its norm belongs
to $W^{1,2}(M;\R)$, and the refined Kato inequality for closed and
coclosed forms \cite[Theorems~3.1(ii) and~6.3(ii)]{CGH} gives 
\[
|d|\omega||^{2}\leq\frac{q}{q+1}|\nabla\omega|^{2}
\]
almost everywhere.
The preceding inequalities, Theorem~\ref{thm:algebraic}, and \eqref{eq:bochner}
imply 
\begin{align*}
0 & <\int_{M}\left(|d|\omega||^{2}+\frac{j(n-j)q}{3n(q+1)}D_{2}|\omega|^{2}\right)\,dV\\
 & \leq\frac{q}{q+1}\int_{M}\bigl(|\nabla\omega|^{2}+\langle Q_{j}(R)\omega,\omega\rangle\bigr)\,dV \\
 & =0,
\end{align*}
a contradiction. Hodge theory proves the real cohomology vanishing,
and the rational homology conclusion follows as in Theorem~\ref{thm:main}. 
\end{proof}

\begin{remark} The coefficient can also be written as 
\[
C_{n}=\frac{\lfloor\frac{n^{2}}{4}\rfloor}{3n}\frac{\lceil\frac{n}{2}\rceil}{\lceil\frac{n}{2}\rceil+1}.
\]
It is smaller than the coefficient $\frac{\lfloor\frac{n^{2}}{4}\rfloor}{3n}$
obtained from the ordinary Kato inequality. Positivity with that coefficient
implies \eqref{eq:spectral-pinching} by the same convex-combination
argument, even when $D_{2}$ changes sign.

If $D_{2}\geq0$ everywhere and is positive at some point, then \eqref{eq:spectral-pinching}
holds. Indeed, the Rayleigh infimum is attained. If it were zero,
a minimizing function would have zero Dirichlet energy and would therefore
be constant. Its potential energy would then force this constant to
vanish, a contradiction. \end{remark}


\begin{thebibliography}{10}
\bibitem{Berger} M. Berger, \emph{Sur les groupes d'holonomie homogène
des variétés à connexion affine et des variétés riemanniennes}, Bull.
Soc. Math. France \textbf{83} (1955), 279-330. \href{https://doi.org/10.24033/bsmf.1464}{doi:10.24033/bsmf.1464}.

\bibitem{Berger1966} M. Berger, \emph{Sur les variétés d'Einstein
compactes}, in Comptes Rendus de la IIIe Réunion du Groupement des
Mathématiciens d'Expression Latine (Namur, 1965), Librairie Universitaire,
Louvain, 1966, pp.~35-55.

\bibitem{Besse} A. L. Besse, \emph{Einstein Manifolds}, Ergebnisse
der Mathematik und ihrer Grenzgebiete (3), vol.~10, Springer-Verlag,
Berlin, 1987. \href{https://doi.org/10.1007/978-3-540-74311-8}{doi:10.1007/978-3-540-74311-8}.

\bibitem{BM} R. G. Bettiol and R. A. E. Mendes, \emph{Sectional curvature
and Weitzenböck formulae}, Indiana Univ. Math. J. \textbf{71} (2022),
no.~3, 1209-1242. \href{https://doi.org/10.1512/iumj.2022.71.8927}{doi:10.1512/iumj.2022.71.8927}.

\bibitem{BC} R. L. Bishop and R. J. Crittenden, \emph{Geometry of
Manifolds}, Pure and Applied Mathematics, vol.~15, Academic Press,
New York-London, 1964.

\bibitem{BS}
S. Brendle and R. M. Schoen,
\emph{Classification of manifolds with weakly $\frac{1}{4}$-pinched curvatures},
Acta Math. \textbf{200} (2008), no.~1, 1-13.
\href{https://doi.org/10.1007/s11511-008-0022-7}{doi:10.1007/s11511-008-0022-7}.

\bibitem{CGH} D. M. J. Calderbank, P. Gauduchon, and M. Herzlich,
\emph{Refined Kato inequalities and conformal weights in Riemannian
geometry}, J. Funct. Anal. \textbf{173} (2000), no.~1, 214-255. \href{https://doi.org/10.1006/jfan.2000.3563}{doi:10.1006/jfan.2000.3563}.

\bibitem{CR} E. Costa and E. Ribeiro, Jr., \emph{Four-dimensional
compact manifolds with nonnegative biorthogonal curvature}, Michigan
Math. J. \textbf{63} (2014), no.~4, 747-761. \href{https://doi.org/10.1307/mmj/1417799224}{doi:10.1307/mmj/1417799224}.

\bibitem{DGM} M. Domínguez-Vázquez, D. González-Álvaro, and L. Mouillé,
\emph{Infinite families of manifolds of positive $k$th-intermediate
Ricci curvature with $k$ small}, Math. Ann. \textbf{386} (2023),
1979-2014. \href{https://doi.org/10.1007/s00208-022-02420-w}{doi:10.1007/s00208-022-02420-w}.

\bibitem{Gallot} S. Gallot, \emph{Estimées de Sobolev quantitatives
sur les variétés riemanniennes et applications}, C. R. Acad. Sci.
Paris Sér. I Math. \textbf{292} (1981), no.~6, 375-377.

\bibitem{GM} S. Gallot and D. Meyer, \emph{Opérateur de courbure
et laplacien des formes différentielles d'une variété riemannienne},
J. Math. Pures Appl. (9) \textbf{54} (1975), no.~3, 259-284.

\bibitem{Ge} J. Ge, \emph{Sharp Weitzenböck and PIC2 estimates from
sectional-scalar curvature pinching}, preprint, 2026. \href{https://arxiv.org/abs/2607.18216}{arXiv:2607.18216}.

\bibitem{Gray} A. Gray, \emph{Compact Kähler manifolds with nonnegative
sectional curvature}, Invent. Math. \textbf{41} (1977), no.~1, 33-43.
\href{https://doi.org/10.1007/BF01390163}{doi:10.1007/BF01390163}.

\bibitem{GX} J.-R. Gu and H.-W. Xu, \emph{The sphere theorems for
manifolds with positive scalar curvature}, J. Differential Geom. \textbf{92}
(2012), no.~3, 507-545. \href{https://arxiv.org/abs/1102.2424}{arXiv:1102.2424}.

\bibitem{GW} L. Guijarro and F. Wilhelm, \emph{Focal radius, rigidity,
and lower curvature bounds}, Proc. Lond. Math. Soc. (3) \textbf{116}
(2018), no.~6, 1519-1552. \href{https://doi.org/10.1112/plms.12113}{doi:10.1112/plms.12113}.

\bibitem{GWConnectivity} L. Guijarro and F. Wilhelm, \emph{A softer
connectivity principle}, Comm. Anal. Geom. \textbf{30} (2022), no.~5,
1093-1119. \href{https://doi.org/10.4310/CAG.2022.v30.n5.a5}{doi:10.4310/CAG.2022.v30.n5.a5}.

\bibitem{KLZ} S. Kong, P. Li, and D. Zhou, \emph{Spectrum of the
Laplacian on quaternionic Kähler manifolds}, J. Differential Geom.
\textbf{78} (2008), no.~2, 295-332. \href{https://doi.org/10.4310/JDG/1203000269}{doi:10.4310/JDG/1203000269}.

\bibitem{Li} X. Li, \emph{Sectional curvature, isotropic curvature,
and Yau's pinching problem}, J. Geom. Anal. \textbf{36} (2026), no.~8,
article~269. \href{https://doi.org/10.1007/s12220-026-02518-6}{doi:10.1007/s12220-026-02518-6}.

\bibitem{LiCOSK} X. Li, \emph{New sphere theorems under curvature
operator of the second kind}, J. Lond. Math. Soc. (2) \textbf{112}
(2025), no.~5, e70356. \href{https://doi.org/10.1112/jlms.70356}{doi:10.1112/jlms.70356}.

\bibitem{LMPUpper} X. Li, E. Matinpour, and P. Petersen, \emph{Homology
sphere theorems from upper bounds for sectional curvature}, manuscript.

\bibitem{LPCOSK} X. Li and P. Petersen, \emph{New curvature operator
of the second kind conditions for rigidity}, preprint, 2026. \href{https://arxiv.org/abs/2610.07687}{arXiv:2610.07687}.

\bibitem{MSW} O. Macia, U. Semmelmann, and G. Weingart, \emph{On
quaternionic bisectional curvature}, Math. Ann. \textbf{391} (2025),
4389-4407. \href{https://doi.org/10.1007/s00208-024-03028-y}{doi:10.1007/s00208-024-03028-y}.

\bibitem{Magnusson} G. {Þ}ór Magnússon, \emph{Integrating holomorphic
sectional curvatures}, preprint, 2023. \href{https://arxiv.org/abs/2307.01332}{arXiv:2307.01332}.

\bibitem{Meyer} D. Meyer, \emph{Sur les variétés riemanniennes à
opérateur de courbure positif}, C. R. Acad. Sci. Paris Sér. A-B \textbf{272}
(1971), A482-A485.

\bibitem{Mouille} L. Mouillé, \emph{Torus actions on manifolds with
positive intermediate Ricci curvature}, J. Lond. Math. Soc. (2) \textbf{106}
(2022), no.~4, 3792-3821. \href{https://doi.org/10.1112/jlms.12674}{doi:10.1112/jlms.12674}.

\bibitem{NPW} J. Nienhaus, P. Petersen, and M. Wink, \emph{Betti
numbers and the curvature operator of the second kind}, J. Lond. Math.
Soc. (2) \textbf{108} (2023), no.~4, 1642-1668. \href{https://doi.org/10.1112/jlms.12790}{doi:10.1112/jlms.12790}.

\bibitem{Nikolayevsky} Y. Nikolayevsky, \emph{Osserman Conjecture
in dimension $n\ne8,16$}, Math. Ann. \textbf{331} (2005), 505-522.
\href{https://doi.org/10.1007/s00208-004-0580-8}{doi:10.1007/s00208-004-0580-8}.
\href{https://arxiv.org/abs/math/0204258}{arXiv:math/0204258}.

\bibitem{NW}
L. Ni and B. Wilking,
\emph{Manifolds with $\frac{1}{4}$-pinched flag curvature},
Geom. Funct. Anal. \textbf{20} (2010), no.~2, 571-591.
\href{https://doi.org/10.1007/s00039-010-0068-5}{doi:10.1007/s00039-010-0068-5}.

\bibitem{Petersen} P. Petersen, \emph{Riemannian Geometry}, third
edition, Graduate Texts in Mathematics, vol.~171, Springer, Cham,
2016. \href{https://doi.org/10.1007/978-3-319-26654-1}{doi:10.1007/978-3-319-26654-1}.

\bibitem{PT} P. Petersen and T. Tao,
\emph{Classification of almost quarter-pinched manifolds},
Proc. Amer. Math. Soc. \textbf{137} (2009), no.~7, 2437-2440.
\href{https://doi.org/10.1090/S0002-9939-09-09802-5}{doi:10.1090/S0002-9939-09-09802-5}.

\bibitem{PW} P. Petersen and M. Wink, \emph{New curvature conditions
for the Bochner technique}, Invent. Math. \textbf{224} (2021), no.~1,
33-54. \href{https://doi.org/10.1007/s00222-020-01003-3}{doi:10.1007/s00222-020-01003-3}.

\bibitem{RW} P. Reiser and D. J. Wraith, \emph{Intermediate Ricci
curvatures and Gromov's Betti number bound}, J. Geom. Anal. \textbf{33}
(2023), no.~12, article~364. \href{https://doi.org/10.1007/s12220-023-01423-6}{doi:10.1007/s12220-023-01423-6}.

\bibitem{SY} R. Schoen and S.-T. Yau, \emph{Lectures on Differential
Geometry}, Conference Proceedings and Lecture Notes in Geometry and
Topology, vol.~I, International Press, Cambridge, MA, 1994.

\bibitem{Shen} Z. Shen, \emph{On complete manifolds of nonnegative
$k$th-Ricci curvature}, Trans. Amer. Math. Soc. \textbf{338} (1993),
no.~1, 289-310. \href{https://doi.org/10.1090/S0002-9947-1993-1112548-6}{doi:10.1090/S0002-9947-1993-1112548-6}.

\bibitem{Simons} J. Simons, \emph{On the transitivity of holonomy
systems}, Ann. of Math. (2) \textbf{76} (1962), no.~2, 213-234. \href{https://doi.org/10.2307/1970273}{doi:10.2307/1970273}.

\bibitem{Wu} H. Wu, \emph{Manifolds of partially positive curvature},
Indiana Univ. Math. J. \textbf{36} (1987), no.~3, 525-548. \href{https://doi.org/10.1512/iumj.1987.36.36029}{doi:10.1512/iumj.1987.36.36029}.

\bibitem{XG} H.-W. Xu and J.-R. Gu, \emph{Rigidity of Einstein manifolds
with positive scalar curvature}, Math. Ann. \textbf{358} (2014), no.~1-2,
169-193. \href{https://doi.org/10.1007/s00208-013-0957-7}{doi:10.1007/s00208-013-0957-7}.

\bibitem{Yau} S.-T. Yau, \emph{Open problems in geometry}, in Differential
Geometry: Partial Differential Equations on Manifolds (Los Angeles,
CA, 1990), Proc. Sympos. Pure Math. \textbf{54}, Part~1, Amer. Math.
Soc., Providence, RI, 1993, pp.~1-28. \href{https://doi.org/10.1090/pspum/054.1/1216573}{doi:10.1090/pspum/054.1/1216573}.

\end{thebibliography}
\end{document}